\documentclass[11pt]{article}
\usepackage{amsfonts,latexsym,amsmath,amscd,geometry}
\usepackage{amssymb}
\usepackage{latexsym}

\newcommand \nc{\newcommand}
\newtheorem{theorem}{Theorem}[section]
\newtheorem{lemma}[theorem]{Lemma}
\newtheorem{corollary}[theorem]{Corollary}

\newtheorem{remark}[theorem]{Remark}
\numberwithin{equation}{section}

\nc{\ba}{\begin{array}}\nc{\ea}{\end{array}}
\nc{\be}{\begin{eqnarray}}\nc{\ee}{\end{eqnarray}}
\nc{\beq}{\begin{equation}}\nc{\eeq}{\end{equation}}
\nc{\bex}{\begin{eqnarray*}}\nc{\eex}{\end{eqnarray*}}
\nc{\btm}{\begin{theorem}} \nc{\etm}{\end{theorem}}
\nc{\blm}{\begin{lemma}} \nc{\elm}{\end{lemma}}
\nc{\R}{\mathbb{R}} \nc{\va}{\varepsilon} \nc{\ls}{\limits}

\def\pf{\noindent{\bf Proof.\quad}}

\newcommand \qed {\hfill $\Box$}

\begin{document}
\title{Harmonic maps in connection of phase transitions with higher dimensional potential wells}
\author{Fanghua Lin\footnote{Courant Institute of Mathematical Sciences, New York University, NY 10012.}
\ and\  Changyou Wang \footnote{Department of Mathematics, Purdue University, West Lafayette, IN 47907.}
\footnote{Both authors are partially supported by NSF.}}
\date{}
\maketitle

\centerline{\it Dedicated to Professor Andy Majda with Deep Admiration}

\begin{abstract}
This is in the sequel of authors' paper \cite{LPW} in which we had set up a program to verify 
rigorously some formal statements associated with the multiple component 
phase transitions with 
higher dimensional wells. The main goal here is to establish a regularity theory for minimizing maps with a 
rather non-standard boundary condition at the sharp interface of the transition. We also present a proof, under 
simplified geometric assumptions, of existence of local smooth gradient flows 
under such constraints on 
interfaces which are in the motion by the mean-curvature. In a forthcoming 
paper, a general theory for such 
gradient flows and its relation to Keller-Rubinstein-Sternberg's work 
\cite{KRS1, KRS2} on the fast reaction, 
slow diffusion and motion by the mean curvature would be addressed.
\end{abstract}

\section {Introduction}
\setcounter{equation}{0}
\setcounter{theorem}{0}
This is a continuation of our previous work Lin-Pan-Wang \cite{LPW} in which we had set up a program to
verify various phenomena associated with multiple components phase transitions with higher dimensional
wells. One of the goals here is to show rigorously the formal asymptotic arguments for the description of 
fast reaction, slow diffusion and sharp interface dynamics using the Ginzburg-Landau approximation as in the 
celebrated papers \cite{KRS1, KRS2} by Keller-Rubinstein-Sternberg. For the leading term of the energy functional 
in the static energy minimization, we showed in \cite{LPW} that the sharp interfaces for these general phase 
transition problem must be area minimizing hypersurfaces with weights. For the energy minimization, each of  
weights must be a constant giving by the length of a so-called minimal connection between a pair of potential 
wells. Therefore for the gradient flow, the dynamic of these sharp interfaces would simply be the motion by 
mean curvature provided that this weight function remains to be a constant that equals the length of a minimal 
connection. The latter leads to a challenging issue of studying energy minimizing maps (phases) and its 
gradient flows that lie in multiple potential wells (submanifolds) of high dimensions and, that each patch of 
such maps (phases) possesses a specific and non-standard boundary condition at corresponding sharp interfaces. The 
phases and their dynamics within each of the potential wells would be derived from the ``slow diffusion" part as 
in \cite{KRS1, KRS2}, and it is hence in the next term of formal asymptotic for the energy of the system. This 
gives a nonlinear coupling between terms of different orders (in formal expansions) of the energy through 
boundary conditions, and it leads us to the study of harmonic maps with these unusual boundary conditions. In 
this paper, we show a boundary regularity theory of minimizing harmonic maps in the above described problems.
We also establish a theorem on the short time existence of classical solutions to the corresponding heat flows.
In a forthcoming work, we will address these dynamical issues in a more general context.

Let us first recall the Cahn-Hilliard energy functional that models the phase transition described by a scalar
function $v$:
$$E_\epsilon(v)=\int_\Omega \big(\epsilon |\nabla v|^2+\frac{1}{\epsilon}W(v)\big)\,dx,$$
where $\Omega\subset\mathbb R^n$ is assumed to be a bounded, smooth domain in $\mathbb R^n$
throughout this paper, $v:\Omega\mapsto\mathbb R$ is the density function, and
$W:\mathbb R\mapsto\mathbb R_+$ is a double-well potential function that has two minima
(zeros) at $\pm 1$. The term $\epsilon |\nabla v|^2$ is the interfacial energy that penalizes the formation
of interface. The asymptotic behavior of minimizers $v_\epsilon$ of $E_\epsilon(\cdot)$ under the constraint
$\int_\Omega v_\epsilon=c$, as $\epsilon\rightarrow 0$, was first studied by Modica and Mortola \cite{MM},
Modica \cite{M}, and Luckhaus and Modica \cite{LM}: they have showed that the separation region between the two stable
phases has $O(\epsilon)$ thickness and the phase transition converges to a minimal hypersurface within the frame work of
De Giorgi's $\Gamma$-convergence theory. There are many important contributions to this problem,
see for examples \cite{FT, KS, MM, M, LM, P1,P2}. 

Rubinstein, Sternberg, and Keller \cite{KRS1, KRS2} introduced the vector-valued system of fast reaction and slow diffusion:
$$\partial_t v_\epsilon=\epsilon\Delta v_\epsilon-\epsilon^{-1}W_v(v_\epsilon) \ {\rm{in}}\ \Omega;
\ \ \frac{\partial v_\epsilon}{\partial\nu}=0 \ {\rm{on}}\ \partial\Omega,$$
where the order paramter $v_\epsilon:\Omega\mapsto\mathbb R^k$ represents the multiple component
phases, and $W:\mathbb R^k\mapsto\mathbb R_+$ vanishes on two disjoint submanifolds in $\mathbb R^k$.
In this case, a front develops in $\Omega$. By the formal WKB analysis on the asymptotic expansion for potential functions
vanishing on two submanifolds, it was found in \cite{KRS1, KRS2} the front moves by its mean curvature, and $v_\epsilon$
approximates the heat flow of harmonic maps away from the front.
Although there have been many studies for the rigorous analysis of such an asymptotics for the scalar case $k=1$,
the corresponding analysis has remained an open problem for $k\ge 2$.

Next we recall the main results of \cite{LPW}. For $k> 1$, let $$N=N^+\cup N^-\subset \mathbb R^k$$
be the union of two disjoint, compact, connected, smooth Riemannian manifolds $N^\pm
\subset\mathbb R^k$ without boundaries.  For $\delta>0$, let
$$N_\delta=\big\{p\in \mathbb R^k:  \ d(p,N)=\inf_{y\in N} |p-y|\le\delta\big\}$$
denote the $\delta$-neighborhood of $N$. It is well known that there exists $\delta_N>0$ such that
$d^2(p, N)\in C^\infty(N_{\delta_N})$. Consider the class of double-well potential functions depending 
only on the distance function from $N$, namely, 
$$F(p)=f(d^2(p,N)),$$
where $f\in C^\infty(\mathbb R_+,\mathbb R_+)$ satisfies the property that there exist
$c_1, c_2, c_3>0$ such that
\begin{equation}\label{f-cond}
\begin{cases}
c_1 t\le f(t)\le c_2 t & \ {\rm{if}}\ 0\le t\le \delta_N^2,\\
f(t)\ge c_3 & \ {\rm{if}}\ 0\le t\le 4\delta_N^2.
\end{cases}
\end{equation}
Consider the family of Cahn-Hiliard functional 
$$E_\epsilon(u)=\int_\Omega (\epsilon^2 |\nabla u|^2+F(u))\,dx, \ u\in H^1(\Omega,\mathbb R^k), \ \epsilon>0,$$
that are singular perturbations of the functional of phase transitions of high dimensional wells:
$$E_0(u)=\int_\Omega F(u)\,dx, \ u\in L^1(\Omega, \mathbb R^k).
$$
For the boundary conditions, we let $\Sigma^\pm\subset\partial\Omega$ be two disjoint, connected,
open subsets of $\partial\Omega$ such that\\
(1) $\partial\Sigma^+=\partial\Sigma^-=\Sigma$ is a connected $(n-2)$-dimensional smooth manifold and\\
(2) $\partial\Omega=\Sigma^+\cup\Sigma^-\cup\Sigma$.

For any small $\eta>0$, let 
$\Sigma^\eta=\big\{x\in\mathbb R^n: \ d(x,\Sigma)<\eta\big\}$ be the $\eta$-neighborhood of $\Sigma$, and
denote $\Sigma_\eta^\pm=\Sigma^\pm\setminus \Sigma^\eta$. Assume that for some $\beta>0$,
$R>0$, $L>0$, and $C>0$, $g_\epsilon: \partial\Omega\mapsto\mathbb R^k$ satisfy:
\begin{itemize}
\item [(1)] $g_\epsilon(\Sigma^\pm_{\epsilon^\beta})\subset N^\pm, g_\epsilon(\partial\Omega)\subset B_R^k,$
and
\begin{equation}\label{ass1}
\int_{\partial\Omega}(\epsilon|\nabla_\tau g_\epsilon|^2+\epsilon^{-1}F(g_\epsilon))\,d\sigma\le L;
\end{equation}
\item [(2)] for any $p^\pm\in N^\pm$, $\exists$ extension maps
$$G_\epsilon^\pm: \Sigma^\pm_{\epsilon^\beta}\times [0,\epsilon^\beta]\mapsto N^\pm$$
such that
$$G_\epsilon^\pm\big|_{\Sigma^\pm_{\epsilon^\beta}\times\{0\}}=g_\epsilon, 
\  G_\epsilon^\pm\big|_{\Sigma^\pm_{\epsilon^\beta}\times\{\epsilon^\beta\}}=p^\pm,$$
\begin{equation}\label{ass2}
\int_{\Sigma^\pm_{\epsilon^\beta}\times [0,\epsilon^\beta]} |\nabla G_\epsilon^\pm|^2\,dx
\le C\Big\{\epsilon^\beta \int_{\Sigma^\pm_{\epsilon^\beta}}|\nabla_\tau g_\epsilon|^2\,dH^{n-1}
+\frac{1}{\epsilon^\beta}\int_{\Sigma^\pm_{\epsilon^\beta}}|g_\epsilon-p^\pm|^2\,dH^{n-1}\Big\},
\end{equation}
where $\nabla_\tau$ denotes the tangential derivative on hypersurfaces in $\mathbb R^n$.
\end{itemize}
Set 
\begin{equation}\label{CH}
{\bf E}(\epsilon)=\min\Big\{\int_\Omega \big(|\nabla u|^2+\frac{1}{\epsilon^2}F(u))\,dx:
\ u\big|_{\partial\Omega}=g_\epsilon\Big\}.
\end{equation}
In \cite{LPW}, we proved 

\noindent{\bf Theorem} A. {\it Assume that $F\in C^\infty(\mathbb R^k)$ satisfies \eqref{f-cond},
$\Gamma\subset\Omega$ is an area-minimizing hypersurface with $\partial\Gamma=\Sigma$
and $g_\epsilon:\partial\Omega\mapsto \mathbb R^k$ satisfies conditions \eqref{ass1}
and \eqref{ass2}. Then
\begin{equation}\label{asym1}
\lim_{\epsilon\rightarrow 0}\epsilon{\bf E}(\epsilon)=c_0^F H^{n-1}(\Gamma),
\end{equation}
where $c_0^F$ is the energy of the minimal connecting orbits between $N^+$ and $N^-$
defined by
\begin{equation}\label{mini_orbit1}
c_0^F=\inf\Big\{c^F(p^+, p^-): \ p^\pm\in N^\pm\Big\},
\end{equation}
and
\begin{equation}\label{mini_orbit2}
c^F(p^+,p^-)
=\inf\Big\{\int_{\mathbb R} (|\xi'(t)|^2+F(\xi))\,dt: \xi\in H^1(\mathbb R,\mathbb R^k), \ \xi(\pm\infty)=p^\pm\Big\}.
\end{equation}
}

Let
$$d_N=\inf\big\{|p^+-p^-|: \ p^\pm\in N^\pm\big\}$$
be the euclidean distance between $N^+$ and $N^-$, and
\begin{equation}\label{mini_set}
\begin{cases}
M^+=\big\{p^+\in N^+: \ \exists\ p^-\in N^-\ s.t.\ |p^+-p^-|=d_N\big\};\\
M^-=\big\{q^-\in N^-: \ \exists\ q^+\in N^+\ s.t.\ |q^+-q^-|=d_N\big\}
\end{cases}
\end{equation}
be the pair of minimal sets in $N^\pm$.

Assume that $g_\epsilon$ is almost optimal near $\Sigma$ in the sense
that its limit $g=\lim_{\epsilon\rightarrow 0}g_\epsilon$
gives the minimal connecting orbits between $N^+$ and $N^-$, 
see \cite{LPW} page 9 for more details.  Then
we also proved in \cite{LPW} the following result.

\noindent{\bf Theorem} B. {\it  Assume $F(p)=f(d^2(p,N))$ satisfies
\eqref{f-cond}, $\Gamma$ is a unique area minimizing hypersurface
with $\partial\Gamma=\Sigma$, which is smooth and strictly stable. 
Assume also that
$${\bf A}=\Big\{v\in H^1(\Omega^\pm, N^\pm):
\ v\big|_{\partial \Omega}=g, \ |v(x^+)-v(x^-)|=d_N\ {\rm{a.e.}}\ x\in\Gamma\big\}\not=\emptyset.$$
Then 
\begin{equation}\label{asymp2}
{\mathbf E}(\epsilon)=\frac{c_0^F}{\epsilon} H^{n-1}(\Gamma)+{\bf D}+o(1),
\end{equation}
where 
\begin{equation}\label{mini_hm}
{\bf D}=\inf\Big\{\int_{\Omega^+} |\nabla v|^2\,dx+\int_{\Omega^-} |\nabla v|^2\,dx: \ v\in{\bf A}\Big\}.
\end{equation}
Furthermore, if $\{u_\epsilon\}$ is a sequence of minimizers of ${\bf E}(\epsilon)$, then
there exists $u\in{\mathbf A}$ attaining the value ${\bf D}$ such that after taking a possible
subsequence, $u_\epsilon$ converges to $u$ in $L^1(\Omega,\mathbb R^k)$.
}

\medskip
The first aim of this paper is to study the boundary regularity of a minimizing harmonic map
$v\in \mathbf{A}$ that attains ${\bf D}$ near the sharp interface $\Gamma$.
In order to achieve it, we make some further assumptions on the minimal sets $M^\pm$.
More precisely, let $M^+\subset N^+$ and $M^{-}\subset N^{-}$
be such that
\begin{itemize}
\item $M^+$ and $M^-$ are connected, $C^1$-manifolds without boundaries, equipped 
with induced metric from $N^+$ and $N^-$ respectively; and 
\item  there exists a $C^1$ differmorphism  $\Phi^+: M^+\mapsto M^-$, whose inverse
map is $\Phi^{-}: M^{-}\mapsto M^+$.

\end{itemize}
Let $\Gamma\subset\Omega$ be a smooth hypersurface with boundary $\Sigma$, i.e. $\partial\Gamma=\Sigma$.
Denote the
two connected components of $\Omega$ separated by $\Gamma$ by $\Omega^\pm$, i.e. $\Omega\setminus\Gamma=\Omega^+\cup\Omega^{-}$,
so that 
$$\partial\Omega^\pm=\Sigma^\pm\cup\Gamma.$$

Let $g:\partial\Omega\to N$ be a given map such that $g\in H^1(\Sigma^\pm, N^\pm)$, and
the two one-side trace values of $g$ on $\Sigma$ satisfy:
\begin{equation}\label{constraint_boundary}
g^\pm(x)(=g(x^\pm))\in H^\frac12(\Sigma, M^\pm), \ \ {\rm{and}}\ \  \Phi^+(g^+(x))=g^{-}(x)\ \ {\rm{a.e.}} \  x\in\Sigma. 
\end{equation}
The minimization problem  seeks
\begin{equation}\label{mini_harmonic}
\inf\Big\{E(u)\ \Big|\
u\in H^1(\Omega^\pm, N^\pm),\  u|_{\partial\Omega}=g, \ u(\Gamma^\pm)\subset M^\pm, \  \Phi^+(u^+(x))=u^{-}(x),
\ {\rm{a.e.}} \ x\in\Gamma\Big\},
\end{equation}
where
$$E(u)=\frac12\int_{\Omega^+}|\nabla u|^2\,dx+\frac12\int_{\Omega^-}|\nabla u|^2\,dx.$$
It is readily seen that if the configuration space
\begin{equation}\label{config_space}
\mathcal A
\equiv\Big\{u\in H^1(\Omega^\pm, N^\pm): \ u|_{\partial\Omega}=g,\ u(\Gamma^\pm)\subset M^\pm, \ 
\ \Phi^+(u^+(x))=u^{-}(x) \ \ {\rm{a.e.}} \ x\in\Gamma\Big\}
\end{equation}
is non-empty,
then there exists at least one energy minimizing map $u\in \mathcal A$, i.e.,
$$E(u)\le E(v),\ \forall v\in\mathcal A.$$

Note that for $n\ge 3$ if, up to a diffeomorphism, $\Omega=B_1\subset \mathbb R^n$, the unit ball,  $\Sigma=\partial B_1\cap \{x_n=0\}$, $\Sigma^\pm=\partial B_1\cap
\mathbb R^n_\pm$, 
$\Gamma= 
B_1\cap \{x_n=0\}$, and $g\in H^1(\Sigma^\pm, N^\pm)$ satisfies (\ref{constraint_boundary}),
then $\mathcal A\not=\emptyset$. In fact, it is not hard to verify that
the homogeneous of degree zero extension 
$\overline{g}(x)=g(\frac{x}{|x|}), x\in B_1$, belongs to $\mathcal A$. In general, we have

\begin{lemma} \label{nonempty}
Assume that $\Pi_1(N^+)=\Pi_1(N^{-})=\{0\}$, $g:\partial\Omega\mapsto N$ satisfies
$g|_{\Sigma^\pm}\in H^1(\Sigma^\pm, N^\pm)$, and the condition \eqref{constraint_boundary} holds.
Then $\mathcal A$ is non-empty.
\end{lemma}
\pf Denote the two one side trace of $g$ on $\Sigma$ by $g^\pm(x)$ for $x\in\Sigma$. Then by (\ref{constraint_boundary})
$g^\pm\in H^\frac12(\Sigma, M^\pm)$. First, we want to extend $g^\pm:\Sigma\mapsto M^\pm$ to maps $G^\pm:\Gamma \mapsto M^\pm$. 
By (\ref{constraint_boundary}), it suffices to construct  an extension map $G^+$ of $g^+$, since
$G^-(x)=\Phi^+(G^+(x))$ for $x\in\Gamma$ will provide an extension of $g^-$. 
Since $M^+$ is connected, i.e. $\Pi_0(M^+)=\{0\}$, Theorem 6.2 of Hardt-Lin \cite{HL2} implies that
for any $1<p<2$, there exists an extension map $G^+\in W^{1,p}(\Gamma, M^+)$ such that
$G^+\big |_{\Sigma}=g^+$ in the trace sense. Now we let $u^+\in H^1(\Omega^+, \mathbb R^k)$ solve
\begin{equation}
\begin{cases}
\Delta u^+=0 \ {\rm{in}}\  \Omega^+\\
\ \ u^+= g \ {\rm{on}} \ \Sigma^+\\
\ \ u^+=G^+ \ {\rm{on}}\ \Gamma.
\end{cases} 
\end{equation}
Since $\Pi_1(N^+)=0$, by applying the extension Lemma 6.1 of \cite{HL2} as in the proof of Theorem 6.2 of \cite{HL2} we conclude that
there exists a map $\widetilde{u}^+\in H^1(\Omega^+, N^+)$ such that
$\widetilde{u}^+-u^+\in H^1_0(\Omega^+, \mathbb R^k)$ and
\begin{eqnarray*}
\int_{\Omega^+}|\nabla\widetilde{u}^+|^2&\le& C\int_{\Omega^+}|\nabla u^+|^2
\le C\big(\|g\|_{H^\frac12(\Sigma^+)}+\|G^+\|_{H^\frac12(\Gamma)}\big)\\
&\le& C\big(\|g\|_{H^\frac12(\Sigma^+)}+\|g^+\|_{H^\frac12(\Sigma)}\big)
\le C\big\|g\big\|_{H^1(\Sigma^+)}.
\end{eqnarray*}
Similarly, we can find an extension map $\widetilde{u}^-\in H^1(\Omega^-, N^-)$ such that
$\widetilde{u}^-=g$ on $\Sigma^{-}$ and $\widetilde{u}^-=G^-$ on $\Gamma$.
Now if we set $\widetilde{u}:\Omega\mapsto N$ by letting $\widetilde{u}(x)=\widetilde{u}^\pm(x)$ for
$x\in\Omega^\pm$, then $\widetilde{u}\in\mathcal A$. This completes the proof.
\qed

\medskip
For a minimizing harmonic map
$u\in\mathcal A$, denote  the set of 
discontinuous points of $u$ in $\Omega^\pm\cup\Gamma$
by $\mathcal S^\pm(u)\subset\Omega^\pm\cup\Gamma$ and define
$$\mathcal S(u)=\mathcal S^+(u)\cup\mathcal S^-(u)$$
as the set of discontinuous points of $u$ in $\Omega$.

It follows from  the interior regularity theory of minimizing harmonic maps by Schoen-Uhlenbeck
\cite{SU1} that $\mathcal S(u)\cap (\Omega\setminus\{\Gamma\})$ has Hausdorff dimension at most $n-3$.

Our first main result concerns the boundary partial regularity at $\Gamma$ for
a minimizing harmonic map $u$ in $\mathcal A$, which is stated as follows.

\begin{theorem}\label{boundary_regularity}
Assume that the boundary value $g\in H^1(\Sigma^\pm, N^\pm)$ satisfies the condition (\ref{constraint_boundary}). If $u\in \mathcal A$ is an energy minimizing harmonic map,
then \\
(i) $\mathcal S(u)\cap\Gamma$ is discrete for $n=3$; and\\
(ii) $\mathcal S(u)\cap\Gamma$ is of Hausdorff dimension at most $(n-3)$ for $n\ge 4$. 
\end{theorem}

The paper is organized as follows. In \S 2, we will give a proof of Theorem \ref{boundary_regularity}. In \S3, we will discuss the corresponding
problem on the heat flow and establish the existence of short time regular solutions. In \S4, we will provide boundary monotonicity inequalities
for both stationary harmonic maps and their corresponding heat flows under the same boundary condition  in Theorem \ref{boundary_regularity},
which may have its own interest and are useful to future studies. 

\section{Proof of Theorem \ref{boundary_regularity}}
\subsection{Euler-Lagrange Equation.} 
In this subsection, we will derive the Euler-Lagrange equation for energy minimizing maps in $\mathcal A$.

Assume that  $u\in\mathcal A$ is an energy minimizing map. For a sufficiently small $\delta>0$, 
let $u(t,\cdot)\in \mathcal A$, $t\in (-\delta,\delta)$, be a family of comparison maps for $u$, i.e. $u(0,\cdot)=u(\cdot)$. For $t\in (-\delta, \delta)$, 
let $u^\pm(t, x)$  denote the two one-sided trace value of $u(t,x)$ for $x\in\Gamma$. Then for $t\in (-\delta,\delta)$, we have
$$u(t, x)=g(x) \ {\rm{for}}\ x\in\Sigma;\  u(t, x)\in N^\pm \ {\rm{for}}\ x\in \Omega^\pm;\  u^\pm(t,x)\in M^\pm\ {\rm{for}}\  x\in\Gamma,$$
and
$$\Phi^+(u^+(t, x))=u^-(t, x)\ \ {\rm{for}}\ H^{n-1}\ {\rm{a.e.}}\ x\in\Gamma.$$
Set $\phi(x)=\frac{d}{dt}\big|_{t=0}u(t, x)$ for $x\in\overline\Omega$. Then we have
\begin{eqnarray*}
0&=&\frac{d}{dt}\big|_{t=0}\big(\frac12\int_{\Omega^+}|\nabla u_t|^2+\frac12\int_{\Omega^-}|\nabla u_t|^2\big)\\
&=&\int_{\Omega^+}\nabla u\cdot\nabla\phi+\int_{\Omega^-}\nabla u\cdot\nabla\phi.
\end{eqnarray*}
For the test function $\phi$, if we denote by $\phi^\pm(x)$ the two one-sided trace value of $\phi$ on $\Gamma$ from $\Omega^\pm$,
then 
$$\phi(x)\in T_{u(x)}N^\pm \ {\rm{for\ a.e.}}\ x\in\Omega^\pm;\ \phi^\pm(x)\in T_{u^\pm(x)} M^\pm
\ \ {\rm{for}}\ H^{n-1} \ {\rm{a.e.}}\ x\in\Gamma, $$
and
$$D\Phi^+(u^+(x))(\phi^+(x))=\phi^-(x) \  \ {\rm{for}}\ H^{n-1} \ {\rm{a.e.}}\ x\in\Gamma.$$
Let $A^\pm$ denote the second fundamental form of $N^\pm$ in $\mathbb R^k$ and denote $u^\pm=u\big|_{\Omega^\pm}$.
Then by integration by parts $u$ satisfies
\begin{equation}\label{HM}
\begin{cases}
-\Delta u^+=A^+(u^+)(\nabla u^+, \nabla u^+)  & {\rm{in}}\ \Omega^+,\\
-\Delta u^-=A^-(u^-)(\nabla u^-, \nabla u^-)  & {\rm{in}}\ \Omega^-,\\
\ \ \ \ \ \ \ u= g  & {\rm{on}}\ \partial\Omega,\\
u^\pm(x)\in M^\pm, \Phi^+(u^+(x))=u^-(x)  & {\rm{on}}\ \Gamma,\\
\big(\frac{\partial u^+}{\partial\nu}\big)^T(x)=(D\Phi^+(u^+(x)))^t\big[\big(\frac{\partial u^-}{\partial\nu}\big)^T(x)\big]
 & {\rm{on}}\ \Gamma.
\end{cases}
\end{equation}
Here $(\cdot)^T(x): T_{u^+(x)}N^+\mapsto T_{u^+(x)}M^+$ (and $(\cdot)^T(x):T_{u^-(x)}N^-\mapsto T_{u^-(x)}M^-$)
denotes the orthogonal projection map for $x\in\Gamma$, and
$$P^t: T_{u^-(x)}M^-\mapsto T_{u^+(x)} M^+ \ ({\rm{or}}\ T_{u^+(x)}M^+\mapsto T_{u^-(x)} M^-)$$
denotes the adjoint of the linear map 
$$P: T_{u^+(x)}M^+\mapsto T_{u^-(x)} M^- \ ({\rm{or }}\ T_{u^-(x)}M^-\mapsto T_{u^+(x)} M^+).$$
It is not hard to see that (2.2)$_5$ can also be written as
$$
\big(\frac{\partial u^-}{\partial\nu}\big)^T(x)=(D\Phi^{-}(u^-(x)))^t\big[\big(\frac{\partial u^+}{\partial\nu}\big)^T(x)\big]
\  {\rm{on}}\ \Gamma.
$$

\subsection{Boundary Monotonicity Inequality}

In order to establish the partial boundary regularity for energy minimizing maps in $\mathcal A$, we need a version of boundary monotonicity inequality.

For $R>0$, denote by $B_R\subset\mathbb R^n$ the ball of radius $R$ and center $0$,
$B_R^\pm=B_R\cap \mathbb R^n_\pm$. 
Since $\Gamma$ is smooth,  there exists $r_0=r_0(\Gamma)>0$
such that for any $x_0\in\Gamma$, $0<r\le r_1:=\frac12\min\{r_0, {\rm{dist}}(x_0,\partial\Omega)\}$,
there exist  $C>0$ and $C^1$-diffeomorphism $\Psi: B_{r}(x_0)=B_r(x_0)\cap\Omega\mapsto B_{r}$ so that
\begin{equation}\label{almost_mono}
\Psi\big(\Omega^\pm\cap B_{r}(x_0)\big)=B_{r}^\pm,
\ \left|D\Psi(x)-\mathbb{I}_n\right|\le C|x-x_0| \ {\rm{for}}\  x\in B_{r}(x_0).
\end{equation}
Here $\mathbb{I}_n$ is the identity matrix of order $n$.
By Fubini's theorem,  $u\in H^1(\partial B_r(x_0)\cap \Omega^\pm, N^\pm)$ for almost all $r\in (0, r_1)$
so that if we define
$$
\widetilde{u}(x)=\begin{cases} u(x)  & x\in\Omega\setminus B_r(x_0),\\
u\big(\Psi^{-1}\big[r\frac{\Psi(x)}{|\Psi(x)|}\big]\big)
& x\in \Omega\cap B_r(x_0),
\end{cases}
$$
then $\widetilde u\in\mathcal A$ is a comparison map for $u$. Thus by the energy minimality,
we have
$$\int_{\Omega^+\cap B_r(x_0)}
|\nabla u|^2 +
\int_{\Omega^-\cap B_r(x_0)}
|\nabla u|^2\le \int_{\Omega^+\cap B_r(x_0)}
|\nabla \widetilde {u}|^2
+\int_{\Omega^-\cap B_r(x_0)}
|\nabla \widetilde{u}|^2.$$
Utilizing (\ref{almost_mono}) and direct calculations, we have that
\begin{eqnarray*}
&&(n-2-Cr) \big(\int_{\Omega^+\cap B_r(x_0)}|\nabla u|^2
+\int_{\Omega^-\cap B_r(x_0)}|\nabla u|^2\big)\\
&\le& r\big(\int_{\Omega^+\cap \partial B_r(x_0)}|\nabla u|^2
+\int_{\Omega^-\cap \partial B_r(x_0)}|\nabla u|^2\big)\\
&-&r\big(\int_{\Omega^+\cap \partial B_r(x_0)}\big|\frac{\partial u}{\partial |x-x_0|}\big|^2
+\int_{\Omega^-\cap \partial B_r(x_0)}\big|\frac{\partial u}{\partial |x-x_0|}\big|^2\big).
\end{eqnarray*}
Therefore,  for any $x_0\in\Gamma$ and $r\in (0, r_1)$, we have that
\begin{eqnarray}\label{almost_mono1}
&&\frac{d}{dr}\big[e^{Cr} r^{2-n}(\int_{\Omega^+\cap B_r(x_0)}|\nabla u|^2
+\int_{\Omega^-\cap B_r(x_0)}|\nabla u|^2)\big]\nonumber\\
&&\ge r^{2-n}\big[\int_{\Omega^+\cap \partial B_r(x_0)}\big|\frac{\partial u}{\partial |x-x_0|}\big|^2
+\int_{\Omega^-\cap \partial B_r(x_0)}\big|\frac{\partial u}{\partial |x-x_0|}\big|^2\big]
\end{eqnarray}
holds, provided $u\in\mathcal{A}$ is an energy minimizing map. In particular, by integrating \eqref{almost_mono1} with respect to $r$, we obtain that
for any $x_0\in\Gamma$ and $0<R_1\le R_2<r_1$,
\begin{eqnarray}\label{almost_mono2}
\displaystyle
&e^{CR_1} R_1^{2-n}\displaystyle\big(\int_{\Omega^+\cap B_{R_1}(x_0)}|\nabla u|^2+\int_{\Omega^-\cap B_{R_1}(x_0)}|\nabla u|^2\big)\nonumber\\
&\ \ +\displaystyle\int_{\Omega\cap (B_{R_2}(x_0)\setminus B_{R_1}(x_0))} |x-x_0|^{2-n}\big |\frac{\partial u}{\partial |x-x_0|}\big |^2\nonumber\\
&\ \le e^{CR_2} R_2^{2-n}\displaystyle\big(\int_{\Omega^+\cap B_{R_2}(x_0)}|\nabla u|^2+\int_{\Omega^-\cap B_{R_2}(x_0)}|\nabla u|^2\big)
\end{eqnarray}
holds for any energy minimizing map $u\in\mathcal{A}$.

\subsection{Boundary Extension Lemma}

A crucial ingredient to prove Theorem \ref{boundary_regularity} is the following boundary extension lemma, 
similar to \cite{HL3} Lemma 3.1. 

\begin{lemma}\label{boundary_extension}
There exist positive constants $\delta$, $q$, and $C$ such that, if
$0<\epsilon<1$, $x_0\in\Gamma$, and
$0<r_0<{\rm{dist}}(x_0,\partial\Omega)$, 
if $\eta^\pm\in H^1(\partial B_{r_0}(x_0)\cap \Omega^\pm, N^\pm)$
satisfies
\begin{eqnarray}\label{small_osc}
&&\int_{\partial B_{r_0}(x_0)\cap\Omega^\pm}
|\nabla_{\rm{tan}}\eta^\pm|^2\,dH^{n-1}
\big[\int_{\partial B_{r_0}(x_0)\cap \Omega^\pm}
|\eta^\pm-p^\pm|^2\,dH^{n-1}
+\int_{\partial B_{r_0}(x_0)\cap\Gamma} |\eta^\pm-p^\pm|^2\,dH^{n-2}\big]\nonumber\\
&&\le\delta^2\epsilon^q
\end{eqnarray}
for some $p^\pm\in\mathbb R^k$, and if $\eta^\pm: {\partial B_{r_0}(x_0)\cap\Gamma} \mapsto M^\pm$  satisfies $$\eta^-(x)=\Phi^+(\eta^+(x))
\ {\rm{for}}\ H^{n-2} \ {\rm{a.e.}}\ x\in \partial B_{r_0}(x_0)\cap\Gamma,$$
then
there exist maps $\omega^\pm\in H^1(B_{r_0}(x_0)\cap\Omega^\pm, N^\pm)$ such that
$\omega^\pm=\eta^\pm$ on $\partial B_{r_0}(x_0)\cap\Omega^\pm$,
and $\omega^\pm: {B_{r_0}(x_0)\cap\Gamma}\mapsto M^\pm$ 
satisfies 
$$\omega^-(x)=\Phi^+(\omega^+(x)) \  {\rm{for}}\  H^{n-1} \ {\rm{a.e.}}\ x\in B_{r_0}(x_0)\cap\Gamma.$$
Furthermore, it holds that
\begin{eqnarray}\label{bdry_extension_estimate}
&&\int_{B_{r_0}(x_0)\cap\Omega^\pm}|\nabla\omega^\pm|^2\,dx
\le \epsilon \int_{\partial B_{r_0}(x_0)\cap\Omega^\pm}
|\nabla_{\rm{tan}}\eta^\pm|^2\,dH^{n-1}\nonumber\\
&&+C\epsilon^{-q}\big[\int_{\partial B_{r_0}(x_0)\cap \Omega^\pm}
|\eta^\pm-p^\pm|^2\,dH^{n-1}
+\int_{\partial B_{r_0}(x_0)\cap\Gamma} |\eta^\pm-p^\pm|^2\,dH^{n-2}\big].
\end{eqnarray}
Here $\nabla_{\rm{tan}}$ denotes the tangential gradient on $\partial B_{r_0}(x_0)$.
\end{lemma}
\pf The proof can be done by suitable modifications of the arguments from \cite{HL2,HL3} and \cite{SU1}. It is based on an induction of the dimension $n$. 
There are two crucial ingredients of the construction: 
\begin{itemize}
\item[(i)]construction in dimension $n=2$; and 
\item[(ii)] homogeneous of degree zero extension for $n\ge 3$. 
\end{itemize}
For simplicity, we will only indicate how to implement these two ingredients in our situation. 
The interested readers can consult with \cite{HL2, HL3, SU1}
for more details.\\
{\it Case} 1: $n=2$ (linear interpolation). Since the problem is invariant under bi-Lipschitz transformations, we may assume that
$x_0=0$, $r_0=1$, $\Omega=B_1$, and $\Gamma=\Gamma_1 (=B_1^2\cap\{x_2=0\})$.  Denote by $S_1^\pm\subset \partial B_1^2$
the half unit circles. Choose $\theta_0^\pm\in S_1^\pm$ so that
$$|\eta^\pm(\theta_0^\pm)-p^\pm|=\inf\Big\{|\eta^\pm(\theta)-p^\pm|: \ \theta\in S_1^\pm\Big\}.$$
Then it is easy to see that
$$\begin{cases}\displaystyle
|\eta^\pm(\theta_0^\pm)-p^\pm|^2\le c\int_{S^\pm_1}|\eta^\pm-p^\pm|^2,\\
\displaystyle\int_{S^\pm_1}|\eta^\pm-\eta^\pm(\theta_0)|^2\le c \int_{S^\pm_1}|\eta^\pm-p^\pm|^2.
\end{cases}
$$
By Sobolev's embedding inequality
$H^1(S_1^\pm)\subset C^\frac12(S_1^\pm)$, we have that
\begin{eqnarray*}
\max_{\theta\in S_1^\pm}\big\{|\eta^\pm(\theta)-\eta^\pm(\theta_0)|^2\big\}
&\le& c\big(\int_{S_1^\pm}|\nabla_{\rm{tan}}\eta^\pm|^2\big)^\frac12
\big(\int_{S_1^\pm}|\eta^\pm-\eta(\theta_0)|^2\big)^\frac12\\
&\le& c\delta\epsilon^{\frac{q}2}.
\end{eqnarray*}
Set 
$$w^+(t,0)=\frac{(1-t)}2 \eta^+(-1,0)+\frac{(1+t)}2 \eta^+(1,0), \ -1\le t\le 1.$$
Then we have
$$\max_{-1\le t\le 1}{\rm{dist}}(w^+(t,0), M^+)\le c|\eta^+(1,0)-\eta^+(-1,0)|\le c\delta\epsilon^{\frac{q}2}.$$
Recall that there exists $\delta_0=\delta_0(M^\pm)>0$ such that for any $0<\delta<\delta_0$, the nearest point
projection maps $\Pi_{M^\pm}: (M^\pm)_\delta\mapsto M^\pm$ and $\Pi_{N^\pm}: (N^\pm)_\delta\mapsto N^\pm$ 
are smooth,  where $(M^\pm)_\delta$ (or $(N^\pm)_\delta$ respectively)
denotes the $\delta$-neighborhood of $M^\pm$  (or $N^\pm$ respectively) in $\mathbb R^k$.
Let $v^+: B_1^+\mapsto\mathbb R^k$ solve
\begin{eqnarray*}
\begin{cases}
\Delta v^+= 0 , &  \ {\rm{in}}\ B_1^+,\\
\ \ v^+= \eta^+, & \ {\rm{on}}\  S_1^+,\\
\ \ v^+=\Pi_{M^+}(w^+), &  \ {\rm{on}}\  \Gamma_1.
\end{cases}
\end{eqnarray*}
Since $\displaystyle\max\big\{{\rm{osc}}_{S_1^+}\eta^+, {\rm{osc}}_{\Gamma_1}\Pi_{M^+}(w^+)\big\}\le C\delta\epsilon^{\frac{q}2}$, 
it follows from the maximum principle  that
$$\max_{x\in B_1^+}{\rm{dist}}(v^+(x), N^+)\le c\delta\epsilon^{\frac{q}2}.$$
Thus we can define
$$\omega^+(x)=\Pi_{N^+}(v^+(x)), \ x\in \overline{B_1^+}.$$
To construct $\omega^-$, first let 
$$w^-(t,0)=\Phi^+(\Pi_{M^+}(w^+(t,0))), \ -1\le t\le 1,$$
so that $w^-(\Gamma_1)\subset M^-$. Let $v^-: B_1^-\mapsto\mathbb R^k$ solve
\begin{eqnarray*}
\begin{cases}
\Delta v^-= 0 , &  \ {\rm{in}}\ B_1^-,\\
\ \ v^-= \eta^-, & \ {\rm{on}}\  S_1^-,\\
\ \ v^-=w^-, &  \ {\rm{on}}\  \Gamma_1.
\end{cases}
\end{eqnarray*}
Then we also have
$$\max_{x\in B_1^-}{\rm{dist}}(v^-(x), N^-)\le c\delta\epsilon^{\frac{q}2},$$
so that we can define
$$\omega^-(x)=\Pi_{N^-}(v^-(x)), \ x\in \overline{B_1^-}.$$
It follows directly from the above construction that 
$\omega^-(x)=\Phi^+(\omega^+(x))$ for $x\in\Gamma_1$, and (\ref{bdry_extension_estimate})
follows from the standard estimate on harmonic functions.

\noindent{\it Case} 2: $n\ge 3$ (homogeneous of degree zero extension). For
$0<\delta<1$, let $B^{\pm,n-1}_{\delta}$ be  $(n-1)$-dimensional half balls of radius 
$\delta>0$, and $C^{\pm, n}_\delta=B^{\pm, n-1}_\delta\times [-\delta,\delta]$ be the 
$n$-dimensional half cylinders of size $\delta$. Let $S_\delta^{\pm, n-2}$ be the $(n-2)$-dimensional half spheres of radius $\delta$ so that
$\partial B_{\delta}^{\pm, n-1}=S_\delta^{\pm,n-2}\cup B_\delta^{n-2}$. 
\begin{lemma} \label{homogeneous} For
$u^\pm\in H^1\big((B_\delta^{\pm, n-1}\times\{\pm\delta\})
\cup (S^{\pm, n-2}_\delta\times [-\delta, \delta]), N^\pm\big)$, if
 $u_1^\pm(x)=u^\pm(x,-\delta)$ and $u_2^\pm(x)=u^\pm(x,\delta)$, $x\in B^{\pm, n-1}_\delta$, 
 satisfies $u_1^\pm, u_2^\pm \in H^1(B_\delta^{\pm, n-1}, N^\pm)$, if $u^\pm(x,t)=u^\pm_0(x)$ for $
(x,t)\in S^{\pm, n-2}_\delta\times [-\delta,\delta]$,
with $u^\pm_0\in H^1(S^{\pm, n-2}_\delta, N^\pm)$, and  if
\begin{equation}\label{bdry_trace}
u^\pm(x,t)\in M^\pm \ {\rm{satisfies}}\ u^-(x,t)=\Phi^+(u^-(x,t)) \ \ {\rm{for}}\ H^{n-2} \ {\rm{a.e.}}\  
x\in B_\delta^{n-2}, \ {\rm{and}}\  t=-\delta, \delta.
\end{equation}
Then there exist extension maps $\overline{u^\pm}\in H^1(C_\delta^{\pm, n}, N^\pm)$ such that
$$\overline{u^\pm}=u^\pm
\ {\rm{ on }}\  (B_\delta^{\pm, n-1}\times\{\pm\delta\})
\cup (S^{\pm, n-2}_\delta\times [-\delta, \delta]),$$ 
\begin{equation}\label{bdry_trace1}
\overline{u^\pm}(x,t)\in M^\pm,
\ \overline{u^-}(x,t)=\Phi^+\big(\overline{u^+}(x,t)\big) \ \ {\rm{for}}\ H^{n-1} \ {\rm{a.e.}}\  
(x,t)\in B_\delta^{n-2}\times [-\delta, \delta],
\end{equation}
and
\begin{eqnarray}\label{extension1}
E\big(\overline{u^\pm}; C_\delta^{\pm, n}\big)
&\le& c\delta\big[E_\delta(u^\pm_1)+E_\delta(u^\pm_2)+\delta E(u^\pm_0)\big],\\
W\big(\overline{u^\pm}; C_\delta^{\pm, n}\big)
&\le& c\delta\big[W_\delta(u^\pm_1)+W_\delta(u^\pm_2)+\delta W(u^\pm_0)\big].
\label{extension2}
\end{eqnarray}
Here $$E_\delta(u_i^\pm)=\int_{B_\delta^{\pm, n-1}}|\nabla u_i^\pm|^2\,dH^{n-1} (i=1,2);
\ E(u_0^\pm)=\int_{S_\delta^{\pm, n-2}}|\nabla_{\rm{tan}}u_0^\pm|^2\,dH^{n-2},$$
and
$$W_\delta(u_i^\pm)=\int_{B_\delta^{\pm, n-1}}|u_i^\pm-p^*|^2\,dH^{n-1} (i=1,2);
\ W(u_0^\pm)=\int_{S_\delta^{\pm, n-2}}|u_0^\pm-p^*|^2\,dH^{n-2}$$ 
for some fixed $p^*\in \mathbb R^L$.
\end{lemma}
{\bf Proof of Lemma 2.3}. By scaling, we may assume $\delta=1$.  There exists  
a bi-Lipschitz homeomorphism $f^\pm: \partial B_1^{\pm, n}\mapsto \partial C_1^{\pm, n}$
such that $\overline{f^\pm}(x)=|x|f^\pm(\frac{x}{|x|}): B_1^{\pm,n}\mapsto C_1^{\pm,n}$ is also
a bi-Lipschitz homeomorphism. Let $\Pi(x)=\frac{x}{|x|}: B_1^n\setminus\{0\}\mapsto \partial B_1^n$ be the radial projection map. 
Define the projection map
$\Pi^\pm: C_1^{\pm, n}\setminus\{0\}\mapsto \partial C_1^{\pm,n}$ by 
$\Pi^\pm=f^\pm\circ\Pi\circ \big(\overline{f^\pm}\big)^{-1}$. Then define
$$\overline{u^\pm}(x)=u^\pm\circ \Pi^\pm(x), \ \ x\in C_1^{\pm,n}.$$
It is easy to see that (\ref{bdry_trace}) implies that $\overline{u^\pm}$ satisfies
the trace condition (\ref{bdry_trace1}) on $\Gamma_1$. It is also easy to see that
$$ E\big(\overline{u^\pm}; C_1^{\pm, n}\big)\le K E\big(\overline{u^\pm}\circ \overline{f^\pm}; B_1^{\pm, n}\big)
\le \frac{K}{n-2} E\big(u^\pm\circ f^\pm; \partial B_1^{\pm, n}\big)
\le C(K) E(u^\pm; \partial C_1^{\pm, n}),
$$
where $K$ is a constant depending on the Lipschitz constants of $f^\pm$ and
$(\overline{f^\pm})^{-1}$. This implies (\ref{extension1}). Similar argument for $W$
also yields (\ref{extension2}).  \qed

\begin{corollary}
There is a constant $c>0$ such that under the same assumptions of Lemma \ref{boundary_extension}, 
if $u\in H^1(\Omega^\pm, N^\pm)\cap\mathcal A$ is energy minimizing among all maps
in $\mathcal A$, and for any $x_0\in\Gamma$ and $0<r_0<{\rm{dist}}(x_0,\partial\Omega)$,
$$r_0^{2-n}\big(\int_{\Omega^+\cap B_{r_0}(x_0)}|\nabla u^+|^2+\int_{\Omega^-\cap B_{r_0}(x_0)}|\nabla u^-|^2\big)
\le c^{-1}\lambda^{-\frac{q}2},$$
then
\begin{eqnarray}\label{bdry_extension_estimate1}
&&(\frac{r_0}2)^{2-n}\big(\int_{\Omega^+\cap B_{\frac{r_0}2}(x_0)}|\nabla u^+|^2
+\int_{\Omega^-\cap B_{\frac{r_0}2}(x_0)}|\nabla u^-|^2\big)\nonumber\\
&&\le \lambda r_0^{2-n}\big(\int_{\Omega^+\cap B_{r_0}(x_0)} |\nabla u^+|^2
+\int_{\Omega^-\cap B_{r_0}(x_0)} |\nabla u^-|^2\big)\nonumber\\
&&+ c\lambda^{-q}r_0^{-n}\big[\int_{\Omega^+\cap B_{r_0}(x_0)}|u^+ - \widehat{u^+}|^2
+\int_{\Omega^-\cap B_{r_0}(x_0)}|u^- - \widehat{u^-}|^2\big]\nonumber\\
&&+ c\lambda^{-q}r_0^{2-n}\big[\int_{\partial B_{r_0}(x_0)\cap\Gamma}|u^+ - \widehat{u^+}|^2\,dH^{n-2}
+\int_{\partial B_{r_0}(x_0)\cap\Gamma}|u^- - \widehat{u^-}|^2\,dH^{n-2}\big],
\end{eqnarray}
where $u^\pm=u|_{\Omega^\pm}$ denotes the restriction of $u$ on $\Omega^\pm$, and 
$\displaystyle\widehat{u^\pm}=\frac{1}{|B_{r_0}(x_0)\cap\Gamma|}\int_{B_{r_0}(x_0)\cap\Gamma} u^\pm\,dH^{n-1}$ 
is the average of the one-side trace of $u^\pm$ in $B_{r_0}(x_0)\cap\Gamma$.
\end{corollary}
\pf  For simplicity, we assume $r_0=1$. Since $u^\pm:\Omega^\pm\to N^\pm$ and $N^\pm$ is compact, it follows 
$$\big|\int_{B_1(x_0)\cap\Gamma} u^\pm\,dH^{n-1}\big|\le C.$$
From the Poincar\'e inequality, we have that
$$\int_{\Omega^\pm\cap B_1(x_0)}|u^\pm-\widehat{u^\pm}|^2\le c\int_{\Omega^\pm\cap B_1(x_0)}|\nabla u^\pm|^2.$$
From the trace estimate and the Poincar\'e inequality, we also have that
$$
\int_{B_1(x_0)\cap\Gamma}|u^\pm-\widehat{u^\pm}|^2\,dH^{n-1}
\le c\big\|u^\pm-\widehat{u^\pm}\big\|_{H^\frac12(B_1(x_0)\cap \Gamma)}^2
\le c\int_{\Omega^\pm\cap B_1(x_0)}|\nabla u^\pm|^2.$$
Applying Fubini's theorem, we can choose $r\in [\frac12,1]$ such that
$$
\int_{\Omega^\pm\cap \partial B_r(x_0)}|\nabla u^\pm|^2\,dH^{n-1}
\le c\int_{\Omega^\pm\cap B_1(x_0)}|\nabla u^\pm|^2
$$
and
\begin{eqnarray*}
&&\big[\int_{\Omega^\pm\cap \partial B_{r}(x_0)}|u^\pm - \widehat{u^\pm}|^2\,dH^{n-1}
+\int_{\partial B_{r}(x_0)\cap\Gamma}|u^\pm - \widehat{u^\pm}|^2\,dH^{n-2}\big]\nonumber\\
&&\le c
\big[\int_{\Omega^\pm\cap B_{1}(x_0)}|u^\pm - \widehat{u^\pm}|^2
+\int_{\Gamma\cap B_{1}(x_0)}|u^\pm - \widehat{u^\pm}|^2\big]\\
&&\le c\int_{\Omega^\pm\cap B_1(x_0)}|\nabla u^\pm|^2.
\end{eqnarray*}
By choosing a sufficiently small $c>0$, we can apply Lemma \ref{boundary_extension}
with $\eta^\pm=u^\pm\big|{\partial B_r(x_0)\cap\Omega^\pm}$ and $p^\pm=\widehat{u^\pm}$
to obtain an extension map $\omega^\pm\in H^1(B_{r}(x_0)\cap\Omega^\pm, N^\pm)$ such that
$\omega^\pm=u^\pm$ on $\partial B_{r}(x_0)\cap\Omega^\pm$,
$\omega^\pm\big|_{B_{r}(x_0)\cap\Gamma}$ has image in $M^\pm$ that
satisfies 
$$\omega^-(x)=\Phi^+(\omega^+(x))
\ {\rm{ for}}\  H^{n-1} \ {\rm{a.e.}}\ x\in B_{r}(x_0)\cap\Gamma,$$
and the estimate (\ref{bdry_extension_estimate}).  If we define $\widetilde u:\Omega\to N$ by
$$\widetilde{u}(x)=
\begin{cases} \omega^\pm(x) & x\in B_r(x_0)\cap\Omega^\pm\\
u(x) & x\in \Omega\setminus B_r(x_0).
\end{cases}
$$
Then $\widetilde u\in\mathcal A$ is a comparison map of $u$. Hence
the energy minimality of $u$ implies that
$$
\int_{\Omega^+\cap B_{r}(x_0)}|\nabla u^+|^2
+\int_{\Omega^-\cap B_{r}(x_0)}|\nabla u^-|^2
\le \int_{\Omega^+\cap B_{r}(x_0)}|\nabla \omega^+|^2
+\int_{\Omega^-\cap B_{r}(x_0)}|\nabla \omega^-|^2,
$$
which, combined with (\ref{bdry_extension_estimate}),  then implies (\ref{bdry_extension_estimate1}). This completes the proof.
\qed

\subsection{Small Energy Regularity}

Another crucial step to prove Theorem \ref{boundary_regularity} is the following energy improvement
property. 
\begin{lemma}\label{energy_improve}
There exist positive constants $\epsilon, C$, and $\theta<1$ such that if $u\in\mathcal A$
is an energy minimizing map that satisfies, for $x_0\in\Gamma$ and some $0<r_0<{\rm{dist}}(x_0,\partial\Omega)$,
\begin{equation}\label{bdry_small_cond}
r_0^{2-n}\big(\int_{\Omega^+\cap B_{r_0}(x_0)} |\nabla u|^2
+\int_{\Omega^-\cap B_{r_0}(x_0)} |\nabla u|^2\big)\leq\epsilon^2,
\end{equation}
then
\begin{eqnarray}\label{bdry_energy_improve}
&&(\theta r_0)^{2-n}\big(\int_{\Omega^+\cap B_{\theta r_0}(x_0)} |\nabla u|^2
+\int_{\Omega^-\cap B_{\theta r_0}(x_0)} |\nabla u|^2\big)\nonumber\\
&\le&\frac12\max\big\{ r_0^{2-n}\big(\int_{\Omega^+\cap B_{r_0}(x_0)} |\nabla u|^2
+\int_{\Omega^-\cap B_{r_0}(x_0)} |\nabla u|^2\big),
\ C{\rm{Lip}}(\Gamma)\big\}.
\end{eqnarray}
\end{lemma}

\medskip
The proof of Lemma \ref{energy_improve} is based on a blowing up argument, similar to \cite{HL3} Theorem 3.3. Before presenting it,  
we need the following regularity estimate on the linear equation, resulting from the blow-up process of the nonlinear harmonic map equation (2.2).

Denote by $B^+_1$ and $B^-_1$ the upper half and lower half unit ball, and
set $\Gamma_1=B_1\cap \{x_n=0\}$. For $a^+\in M^+$, let $a^-=\Phi^+(a^+)\in M^-$.
Let ${\rm{Tan}}({a^\pm}, M^\pm)$ denote the tangent space of $M^\pm$ at $a^\pm$, and
${\rm{Nor}}({a^\pm}, M^\pm)$ denote the normal space of $M^\pm\subset N^\pm$ at $a^\pm$, i.e.,
$${\rm{Tan}}({a^\pm}, M^\pm) \oplus {\rm{Nor}}({a^\pm}, M^\pm)
={\rm{Tan}}({a^\pm}, N^\pm).$$
For any vector $v_\pm\in {\rm{Tan}}({a^\pm}, N^\pm)$, we decompose it as
$$v_\pm =v_\pm^t+v_\pm^n,$$
where $v_\pm^t$ denotes the orthogonal projection of $v_\pm$ into ${\rm{Tan}}({a^\pm}, M^\pm)$, and $v_\pm^n$ denotes the orthogonal projection of $v_\pm$ into ${\rm{Nor}}({a^\pm}, M^\pm)$.

\begin{lemma} \label{linear_regularity} Suppose that $v_\pm\in H^1(B_1^\pm, {\rm{Tan}}(a^\pm, N^\pm))$ are two harmonic functions, with traces 
$v_\pm\big|_{\Gamma_1}\in H^\frac12(\Gamma_1, 
{\rm{Tan}}(a^\pm, M^\pm))$, satisfying
\begin{equation}\label{trace_match}
\begin{cases} v_-=D\Phi^+(a^+)(v_+) &\ {\rm{on}}\ \Gamma_1,\\
\big(\frac{\partial v_+}{\partial x_n}\big)^T= (D\Phi^+(a^+))^t\big(\frac{\partial v_-}{\partial x_n}\big)^T 
&\ {\rm{on}}\ \Gamma_1.
\end{cases}
\end{equation}
Then $v_\pm\in C^\infty\big(B_\frac12^\pm\cup \Gamma_\frac12\big)$, and
for any $l\ge 1$, it holds
\begin{equation}\label{trace_estimate1}
\big\|v_\pm\big\|_{C^l\big(B_\frac12^\pm\cup\Gamma_\frac12\big)}
\le C\big(l, \big\|\Phi^+\big\|_{C^1(M^+)}, \big\|v_\pm\big\|_{H^1(B_1^\pm)}\big).
\end{equation}
\end{lemma}
\pf Since $a^\pm \in M^\pm$, we can decompose
$v_\pm=v_\pm^t+v_\pm^n$ so that
\begin{equation} \label{harm_1}
\Delta v_\pm^t=0 \ \ \ {\rm{in}}\ B_1^\pm,
\end{equation}
and
\begin{equation}\label{harm_2}
\Delta v_\pm^n=0 \ \ \ {\rm{in}}\ B_1^\pm.
\end{equation}
Since $v_\pm(x)\in {\rm{Tan}}(a^\pm, M^\pm)$
for $H^{n-1}$ a.e. $x\in\Gamma_1$, we have that 
\begin{equation}\label{trace_cond}
v_\pm^n =0 \ {\rm{on}} \ \Gamma_1.
\end{equation}
It is readily seen that by (\ref{harm_2}) and (\ref{trace_cond}), 
$v_\pm^n \in C^\infty\big(B_\frac12^\pm\cup\Gamma_\frac12\big)$, and
for any $l\ge 1$
\begin{equation}\label{trace_estimate2}
\big\|v_\pm^n\big\|_{C^l\big(B_\frac12^\pm\cup\Gamma_\frac12\big)}
\le C\big(l, \big\|v_\pm^n\big\|_{H^1(B_1^\pm)}\big).
\end{equation}
To show regularity of $v_\pm^t$, we denote 
$P=D\Phi^+(a^+)$ and proceed as follows. 
Define $\widetilde{v_{-}}:B_1^{+}
\mapsto T_{a^+}N^+$ be an even extension $v_{-}$, i.e.,
$$\widetilde{v_{-}}(x',x_n)=v_{-}(x', -x_n),  \ \ (x', x_n)\in B_1^+.$$
Then it is easy to see that 
\begin{equation}\label{odd_harm}
\begin{cases}
\Delta \big(\widetilde{v_{-}}^t-P(v_{+}^t)\big)=0 & \ {\rm{in}}\ B_1^+,\\
\widetilde{v_{-}}^t-P(v_{+}^t)=0 & \ {\rm{on}}\ \Gamma_1,
\end{cases}
\end{equation}
and
\begin{equation}\label{even_harm}
\begin{cases}
\Delta \big(v_{+}^t-P^t(\widetilde{v_{-}}^t)\big)=0 & \ {\rm{in}}\ B_1^+,\\
\frac{\partial}{\partial x_n}\big(v_{+}^t-P^t(\widetilde{v_{-}}^t)\big)=0 & \ {\rm{on}}\ \Gamma_1.
\end{cases}
\end{equation}
From the standard theory of harmonic functions,
we see that (\ref{odd_harm}) and (\ref{even_harm}) imply
$$\big(\widetilde{v_{-}}^t-P(v_{+}^t)\big),\  \big(v_{+}^t+P^t(\widetilde{v_{-}}^t)\big)
\in C^\infty(B_\frac12^-\cup\Gamma_\frac12),$$
and it holds that, for any $l\ge 1$,
\begin{equation}\label{trace_estimate3}
\big\|\widetilde{v_{-}}^t-P(v_{+}^t)\big\|_{C^l\big(B_\frac12^-\cup\Gamma_\frac12\big)}
+\big\|v_{+}^t+P^t(\widetilde{v_{-}}^t)\big\|_{C^l\big(B_\frac12^-\cup\Gamma_\frac12\big)}
\le C\big(l, \big\|v_{\pm}\big\|_{H^1(B_1^\pm)}\big).
\end{equation}

If $PP^t=I_k$, i.e. $P\in O(k)$ is an orthogonal matrix, then we have
$$|\widetilde{v_{-}}^t-P(v_{+}^t)|=|P(P^t(\widetilde{v_{-}}^t)-v_{+}^t)|=|P^t(\widetilde{v_{-}}^t)-v_{+}^t|.$$
This and \eqref{trace_estimate3} easily yield \eqref{trace_estimate1}.

If $PP^t\not=I_k$, then $P^{-1}\not=P^t$ and we can also see easily that \eqref{trace_estimate1}
follows from \eqref{trace_estimate3}. This completes the proof. \qed

\bigskip
\noindent{\bf Proof of Lemma 2.5}.  The proof follows from a blow-up argument, Lemma \ref{linear_regularity},
and the boundary extension Lemma 2.2. Here we only sketch the argument. 

For simplicity, assume that $x_0=0$, $r_0=1$, $\Omega=B_1$, and $\Gamma=\Gamma_1$ so that ${\rm{Lip}}(\Gamma)=0$. Suppose that the conclusion were false.
Then for any $\theta\in (0,1)$, there would exist $\epsilon_i\rightarrow 0$ and a sequence of minimizing harmonic
maps $u_i\in \mathcal A$ that satisfy
\begin{equation}\label{small_energy_seq}
\int_{B_1^+}|\nabla u_i^+|^2+\int_{B_1^-}|\nabla u_i^-|^2=\epsilon_i^2,
\end{equation}
and
\begin{equation}\label{no_decay}
\theta^{2-n}\big(\int_{B_\theta^+}|\nabla u_i^+|^2
+\int_{B_\theta^-}|\nabla u_i^-|^2\big)>\frac12 \epsilon_i^2.
\end{equation}
Let $\displaystyle\overline{u_i^\pm}=\frac{1}{|\Gamma_1|}\int_{\Gamma_1} u_i^\pm$ denote the average of
the two one-sided traces of $u_i$ on $\Gamma_1$. By the Poincar\'e inequality on $\Gamma_1$
and $H^1$ trace theory, we have
\begin{eqnarray*}
{\rm{dist}}\big(\overline{u_i^+}, M^+\big)^2
\le \frac{1}{|\Gamma_1|}\int_{\Gamma_1}\big|u_i^+-\overline{u_i^+}\big|^2\,dH^{n-1}
\le c\big\|\nabla u_i^+\big\|_{L^2(B_1^+)}^2\le c\epsilon_i^2.
\end{eqnarray*}
Therefore for $i$ sufficiently large there is a unique nearest point $a_i^+=\Pi_{M^+}\big(\overline{u_i^+}\big) \in M^+$ such that
$$\big|\overline{u_i^+}-a_i^+\big|={\rm{dist}}\big(\overline{u_i^+}, M^+\big).$$
Since $u_i^-=\Phi^+(u_i^+)$ on $\Gamma_1$, it is readily seen that 
$a_i^-\equiv \Phi^+(a_i^+)\in M^-$ satisfies
\begin{eqnarray*}\big|\overline{u_i^-}-a_i^-\big|^2
&=&\big|\frac{1}{\Gamma_1}\int_{\Gamma_1}\Phi^+(u_i^+)-\Phi^+(a_i^+)\big|^2
\le c{\rm{Lip}}^2(\Phi^+)\int_{\Gamma_1}|u_i^+-a_i^+|^2\\
&\le& c\int_{\Gamma_1}\big|u_i^+-\overline{u_i^+}\big|^2+\big|\overline{u_i^+}-a_i^+\big|^2\le c\epsilon_i^2.
\end{eqnarray*}
Now we define the corresponding blow-up sequence 
$v_i:B_1\to\mathbb R^k$ by letting
$$
v_i(x):=\begin{cases}v_i^+(x)=\frac{u_i^+(x)-a_i^+}{\epsilon_i}, & x\in B_1^+,\\
v_i^-(x)=\frac{u_i^-(x)-a_i^-}{\epsilon_i}, & x\in B_1^-.
\end{cases}
$$
It is easy to see that 
\begin{equation}\label{H1_bound_vi}
\int_{B_1^+}|\nabla v_i^+|^2+\int_{B_1^-}|\nabla v_i^-|^2=1,
\end{equation}
and
\begin{equation}\label{no_H1_decay}
\theta^{2-n}\big(\int_{B_\theta^+}|\nabla v_i^+|^2+\int_{B_\theta^-}|\nabla v_i^-|^2\big)>\frac12.
\end{equation}
By (\ref{H1_bound_vi}) and the $H^1$-trace theory, we have
$$\displaystyle \|v_i^\pm\|_{H^1(B_1^\pm)}\le c.$$
Hence, after taking a subsequence, there exists $v:B_1\to\mathbb R^k$,
with $v_\pm (=v|_{B_1^\pm}) \in H^1(B_1^\pm,\mathbb R^k)$, such that $v_i^\pm$ converge
to $v_\pm$ weakly in $H^1(B_1^\pm,\mathbb R^k)$. In particular, by (\ref{H1_bound_vi}), we have
\begin{equation}\label{H1_bound_v}
\int_{B_1^+}|\nabla v_+|^2+\int_{B_1^-}|\nabla v_-|^2\le 1.
\end{equation}
Again passing to a subsequence, we assume that
$$\lim_{i\rightarrow\infty} a_i^+=a^+\in M^+
\ {\rm{and}}\ \lim_{i\rightarrow\infty} a_i^-=a^-=\Phi^+(a^+)\in M^-.
$$
It is not hard to verify that $v_{+}(x)\in T_{a^+} N^+$ for a.e. $x\in B_1^+$, and
$v_{-}(x)\in T_{a^-} N^-$ for a.e. $x\in B_1^{-}$.
Since $u_i^\pm(x)\in M^\pm$ for $H^{n-1}$ a.e. $x\in\Gamma_1$, it is also not hard to see that
\begin{equation}\label{trace_cond_v}
v_\pm(x)\in {\rm{Tan}}(a^\pm, M^\pm)
\ {\rm{and}} \ \ v_-(x)=D\Phi^+(a^+)(v_+(x)) 
\  H^{n-1}\ {\rm{a.e.}} \ x\in\Gamma_1.
\end{equation}
Since $v_i^\pm$ satisfies
$$-\Delta v_i^\pm =\epsilon_iA^\pm(u_i^\pm)(\nabla v_i^\pm, \nabla v_i^\pm)
\ \ {\rm{in}} \ \ B_1^\pm,$$
and 
$$\int_{B_1^\pm} |A^\pm(u_i^\pm)(\nabla v_i^\pm, \nabla v_i^\pm)|\le
c\int_{B_1^\pm} |\nabla v_i^\pm|^2\le c,$$
we have, after taking $i$ to infinity, that
\begin{equation}\label{harm_v}
-\Delta v_\pm=0 \ \ {\rm{in}} \ \ B_1^\pm.
\end{equation}
Since $v_i^\pm$ also satisfies the trace condition
$$
\big(\epsilon_i\frac{\partial v_i^{+}}{\partial \nu}\big)^T
=(D\Phi^+(\epsilon_i v_i^+(x)+a_i^+))^t\big(\epsilon_i \frac{\partial v_i^-}{\partial \nu}\big)^T \ \ {\rm{on}}\ \ \Gamma_1,
$$
we obtain, after taking $i$ to infinity,  that 
\begin{equation}\label{trace_cond1_v}
\big(\frac{\partial v_{+}}{\partial \nu}\big)^T
=(D\Phi^+(a^+))^t\big(\frac{\partial v_{-}}{\partial \nu}\big)^T
\ \ {\rm{on}}\ \ \Gamma_1.
\end{equation}
Here $(\cdot)^T: T_{a^\pm}N^\pm \mapsto T_{a^\pm}M^\pm$ is  the orthogonal projection map.
Moreover, we claim 
\begin{equation}\label{zero_average_trace}
\int_{\Gamma_1}v_\pm\,dH^{n-1}=0.
\end{equation}
Set $\displaystyle w_i^+=\frac{\overline{u_i^+}-a_i^+}{\epsilon_i}$ 
and $\displaystyle w^+=\lim_{i\rightarrow\infty} w_i^+$. Then we have that $w^+\in {\rm{Nor}}(a^+, M^+)$. Hence for $H^{n-1}$ a.e. $x\in\Gamma_1$, 
$$\lim_{i\rightarrow \infty}\frac{u_i^+(x)-a_i^+}{\epsilon_i}\cdot w_i^+=0,$$
since $\frac{u_i^+(x)-a_i^+}{\epsilon_i}$ converges to a vector in ${\rm{Tan}}(a^+, M^+)$.
Thus
$$|w^+|^2=\lim_{i\rightarrow\infty}|w_i^+|^2=\frac{1}{|\Gamma_1|}\lim_{i\rightarrow \infty}\int_{\Gamma_1}
\frac{u_i^+(x)-a_i^+}{\epsilon_i}\cdot w_i^+\,dH^{n-1}=0.$$
This implies
$$\frac{1}{|\Gamma_1|}\int_{\Gamma_1}v_+\,dH^{n-1}
=\frac{1}{|\Gamma_1|}\lim_{i\rightarrow\infty}\int_{\Gamma_1}\frac{u_i^+-a_i^+}{\epsilon_i}\,dH^{n-1}
=\lim_{i\rightarrow\infty}\frac{\overline{u_i^+}-a_i^+}{\epsilon_i}=w^+=0.$$
To see $\displaystyle\int_{\Gamma_1} v_-\,dH^{n-1}=0$, observe that
\begin{eqnarray*}\frac{\overline{u_i^-}-a_i^-}{\epsilon_i}
&=&D\Phi^+(a_i^+)\big(\frac{\overline{u_i^+}-a_i^+}{\epsilon_i}\big)
+o(1)\int_{\Gamma_1}|\frac{u_i^+-a_i^+}{\epsilon_i}|\,dH^{n-1}\\
&=&D\Phi^+(a_i^+)(w_i^+)
+o(1)\|v_i^+\|_{L^1(\Gamma_1)}
\end{eqnarray*}
so that
$$\int_{\Gamma_1} v_-\,dH^{n-1}=\frac{1}{|\Gamma_1|}\lim_{i\rightarrow\infty}
\frac{\overline{u_i^-}-a_i^-}{\epsilon_i}=\frac{1}{|\Gamma_1|}
D\Phi^+(a^+)(w^+)=0.$$

By (\ref{trace_cond_v}), (\ref{harm_v}), and (\ref{trace_cond1_v}), we can apply Lemma \ref{linear_regularity} to conclude that 
$v_\pm\in C^\infty\big(\overline{B_\frac12^\pm}\big)$.
Moreover,  by (\ref{H1_bound_v}) and (\ref{zero_average_trace})
we have that for any $0<\theta<1$,
\begin{eqnarray}
&&\theta^{-n}\big(\int_{B_\theta^+}|v_+-(v_+)_\theta|^2
+\int_{B_\theta^-}|v_{-}-(v_{-})_\theta|^2\big)\nonumber\\
&&\le  c\theta^2
\big(\int_{B_1^+}|\nabla v_+|^2
+\int_{B_1^-}|\nabla v_{-}|^2\big)\le c\theta^2,\label{harm_v_estimate1}
\end{eqnarray}
where $\displaystyle
(v_\pm)_\theta=\frac{1}{|\Gamma_\theta|}\int_{\Gamma_\theta} v_\pm \,dH^{n-1}.$
By the Poincar\'e inequality and the trace theory we also have
\begin{equation}\label{harm_v_estimate2}
\theta^{1-n}\big(\int_{\Gamma_\theta}|v_+-(v_+)_\theta|^2
+\int_{\Gamma_\theta}|v_{-}-(v_{-})_\theta|^2\big)
\le c\theta^2\big(\int_{B_1^+}|\nabla v_+|^2
+\int_{B_1^-}|\nabla v_{-}|^2\big)\le c\theta^2 .
\end{equation}
Since $v_i^\pm\rightarrow v_\pm$ in $L^2(B_1^\pm)$
and $L^2(\Gamma_1)$, it follows from (\ref{harm_v_estimate1}) and (\ref{harm_v_estimate2})
that for $i$ sufficiently large 
\begin{eqnarray}
&&\theta^{-n}\big(\int_{B_\theta^+}|u_i^+-(u_i^+)_\theta|^2
+\int_{B_\theta^-}|u_i^--(u_i^-)_\theta|^2\big)
+\theta^{1-n}\int_{\Gamma_\theta}(|u_i^+-(u_i^+)_\theta|^2
+|u_i^--(u_i^-)_\theta|^2)\nonumber\\
&&\le c\theta^2\epsilon_i^2.\label{l2_decay_ui}
\end{eqnarray}
Combining (\ref{bdry_extension_estimate1}) with (\ref{l2_decay_ui}). we can
repeat the argument of \cite{HL2} to get a desired contradiction.

\medskip
\noindent{\bf Proof of Theorem 1.2}. It is well-known that iterations of Lemma \ref{energy_improve}, combined with the interior $\epsilon$-regularity,
implies that there exist $\epsilon_0>0$ and $\alpha_0\in (0,1)$ such that if for $x_0\in\Gamma$, there exists $r_0>0$ such that
$$r_0^{2-n}\big(\int_{\Omega^+\cap B_{r_0}(x_0)} |\nabla u|^2+\int_{\Omega^-\cap B_{r_0}(x_0)} |\nabla u|^2\big)\le \epsilon_0^2,$$
then $u\in C^{\alpha_0}(\overline{\Omega^\pm}\cap B_{\frac{r_0}2}(x_0), N^\pm)$\footnote{Higher order regularity of $u$, e.g, $u\in C^{l,\alpha}(\overline{\Omega^\pm}\cap B_{\frac{r_0}2}(x_0))$,
can be shown, provided that the map $\Phi^+:M^+\to M^{-}$ is assumed to be $C^{l+1,\alpha}$ for some $l\ge 1$ and $0<\alpha<1$.}. 
It follows from this property that  the set $\mathcal{S}(u)$ of discontinuity for $u$ in $\Omega^\pm\cup\Gamma$ can be shown to have
$H^{n-2}(\mathcal{S}(u))=0$. 
It follows from \cite{SU1} that the Hausdorff dimension of $\mathcal{S}(u)$, ${\rm{dim}}_H(\mathcal{S}(u)\cap(\Omega^+\cup\Omega^{-}))\le n-3$ for $n\ge 3$. 
Employing the boundary extension lemma \ref{boundary_extension} and Federer's dimension reduction argument, we can proceed, similar to \cite{SU2} and \cite{HL3}, to  conclude that
${\rm{dim}}_H(\mathcal{S}(u)\cap\Gamma)\le n-3$ for $n\ge 3$, and $\mathcal{S}(u)$ is discrete when $n=3$. This completes the proof. \qed

\section{On the local existence of regular solutions to heat flow}
In this section, we will consider the  gradient flow associated with the minimization problem \eqref{mini_harmonic}, or, equivalently, the parabolic version
of the harmonic map equation \eqref{HM}.  Under some further assumptions on $M^\pm$ and $\Gamma$, to be specified below,
we will establish the local existence of regular solutions of the heat flow under the initial and corresponding boundary conditions.

Before describing the corresponding heat flow problem, we first need to introduce some notations.
For a given $T>0$, let $\{\Gamma(t): t\in [0,T]\}$ be a smooth family of smooth hypersurfaces,
with $\Gamma(0)=\Gamma$, such that  
$$\partial \Gamma(t)=\partial\Gamma=\Sigma, \ \forall 0\le t\le T.$$
For $t\in [0,T]$, decompose $\Omega\setminus\Gamma(t)$ into the disjoint union of
two simply connected components $\Omega^+(t)$ and $\Omega^-(t)$, i.e.,
$$\Omega\setminus\Gamma(t)=\Omega^+(t)\cup \Omega^-(t), \ t\in [0,T].$$
Denote $\Omega^\pm=\Omega^\pm(0)$, and write 
$$\Omega\setminus\Gamma=\Omega^+\cup\Omega^-,\ \ \partial\Omega\setminus\Sigma=\Sigma^+ \cup\Sigma^{-}$$
so that $\partial\Omega^\pm=\Gamma\cup\Sigma^\pm.$
Set 
$$Q_T=\big\{(x,t): \ x\in\Omega, \ 0<t\le T\big\}, \ \partial_pQ_T=(\Omega\times \{0\})\cup (\partial\Omega\times (0,T]),$$
and
$$\Gamma_T=\big\{(x,t): \ x\in \Gamma(t), \ 0<t\le T\big\}, \ Q_T^\pm=\big\{(x,t): \ x\in \Omega^\pm(t), \ 0<t\le T\big\}.$$
The harmonic heat flow problem corresponding to \eqref{HM} can be formulated as follows. We are looking for 
$u^\pm:Q_T^\pm\mapsto N^\pm$, with $u^\pm(x,t)\in M^\pm$ for $(x,t)\in \Gamma_T$,  that solves
\begin{equation}\label{HFHM}
\begin{cases}
\partial_tu^+-\Delta u^+=A^+(u^+)(\nabla u^+, \nabla u^+) & \ {\rm{in}}\ Q_T^+,\\
\partial_tu^--\Delta u^-=A^-(u^-)(\nabla u^-, \nabla u^-) & \ {\rm{in}}\ Q_T^-,\\
u^\pm(x, t)= g^\pm(x) & \ (x,t)\in\Sigma^\pm\times [0, T],\\
u=u_0^\pm & \ {\rm{on}}\ \Omega^\pm\times\{0\},\\
\Phi^+(u^+)=u^-& \ {\rm{on}}\ \Gamma_T,\\
\big(\frac{\partial u^+}{\partial\nu}\big)^T=(D\Phi^+(u^+))^t\big[\big(\frac{\partial u^-}{\partial\nu}\big)^T\big] &  \ {\rm{on}}\ \Gamma_T.
\end{cases}
\end{equation}
Here $u_0^\pm:\Omega^\pm\mapsto N^\pm$, with $u_0^\pm(x)\in M^\pm$ satisfying $u_0^-(x)=\Phi^+(u_0^+(x))$
for $x\in\Gamma$, and $g^\pm=u_0^\pm\big|_{\Sigma^\pm}$ are given initial and boundary values.

In order to establish the short time existence of regular solutions to \eqref{HFHM}, 
we need to set up the problem appropriately by specifying
the assumptions $({\bf A})$, $({\bf B})$, and $({\bf C})$ on $N^\pm$
and $M^\pm$:

\begin{itemize}
\item[({\bf A})] The target Riemannian manifolds $(N^\pm, h^\pm)$ have the same dimension ${\rm{dim}}(N^\pm)=k+m$. 
For, otherwise,
if $k_1={\rm{dim}}(N^+)<k_2={\rm{dim}}(N^-)$, then we can replace $(N^+, h^+)$ by
$$\big(\widehat{N^+}=N^+\times \mathbb S^{k_2-k_1}, \ \widehat{h^+}=h^+\oplus h_{\rm{can}}\big),$$
where $h_{\rm{can}}$ denotes the standard metric on $\mathbb S^{k_2-k_1}$. Notice that ${\rm{dim}}(\widehat{N^+})=k_2$.
Moreover, for any map $u:\Omega^+(t)\times [0,T]\to N^+$,
if we define $\tilde{u}(x,t)=(u(x,t),e): \Omega^+(t)\times [0,T]\to \widehat{N^+}$, where $e\in \mathbb S^{k_2-k_1}$,
then we can show that if $u$ is a solution to the heat flow of harmonic maps to $N^+$, then 
$\tilde{u}$ is also a solution to the heat flow of harmonic maps to $\widehat{N^+}$. This follows from
the chain rule and the fact that $(N^+,h^+)$ is a totally geodesic sub-manifold of $(\widehat{N^+}, \widehat{h^+})$.

\item[({\bf B})] The manifolds  $M^\pm\subset N^\pm$ are two $k$-dimensional compact smooth sub-manifolds,
with $\partial M^\pm=\emptyset$, such that
there exists a smooth  diffeomorphism $\Phi^+: M^+\mapsto M^-$, whose inverse is denoted by $\Phi^-:M^-\mapsto M^+$. 
Moreover, there exists $r_0=r_0(M^+)>0$ such that for any $p^+\in M^+$, $\Phi^+$ can be extended into a smooth
diffeomorphism, still denoted as itself, 
$$\Phi^+: B_{r_0}^{N^+}(p^+)=\{p\in N^+:  {\rm{d}}_{N^+}(p, p^+)<r_0\}\mapsto 
B_{r_0}^{N^-}(p^-)=\{p\in N^-:  {\rm{d}}_{N^-}(p, p^-)<r_0\},$$
whose inverse is also denoted by $\Phi^-$. 

\item[({\bf C})] There exists a $0<r_1=r_1(N^+)\le r_0(M^+)$ such that for any $p^+\in N^+$, 
there exists a local parametrization of $B_{r_1}^{N^+}(p^+)$ by $(B_1^k\times B_1^{m}, \phi^+)$, i.e.,
$$U=(U^1, U^2)=\big((u_1,\cdots, u_k), (u_{k+1},\cdots, u_{k+m})\big)\in B_1^k\times B_1^m$$
provides a local  representation of $B_{r_1}^{N^+}(p^+)$ via the 
diffeomorphism $\phi^+:B_1^k\times B_1^m\mapsto B_{r_1}^{N^+}(p^+)$.
We may assume that $U(p^+)=(0,0)$, and if $p^+\in M^+$ then
$$U\big(M^+\cap B_{r_1}^{N^+}(p^+)\big)\equiv\big\{U=(U^1, U^2)\in B_1^k\times B_1^m: \ U^2=0\big\},$$ 
and the Riemannian metric $h^+$ on $B_{r_1}^{N^+}(p^+)$ can be expressed by
$$h^+(U)=\sum_{i,j=1}^{k+m} h_{ij}^+(U)du_i\otimes du_j,\ \forall U\in B_1^k\times B_1^m,$$
and the induced metric of  $h^+$ on $M^+\cap B_{r_1}^{N^+}(p^+)$ is given by
$$\displaystyle h^+(U^1,0)=\sum_{i,j=1}^{k} h_{ij}^+(U^1,0)du_i\otimes du_j, \ \forall U^1\in B_1^k.$$

\end{itemize}

It is readily seen that for $p^+\in M^+$ and $p^-=\Phi^+(p^+)$, through the diffeomorphism 
$\Phi^+:B_{r_0}^{N^+}(p^+)\mapsto B_{r_0}^{N^-}(p^-)$, $U=(U^1,U^2)\in B_1^k\times B_1^m$ 
provides a local parametrization of $B_{r_1}^{N^-}(p^-)$
through the diffeomorphism $\phi^-:=\Phi^+(\phi^+): B_1^k\times B_1^m\mapsto B_{r_1}^{N^-}(p^-)$. 
In particular, $U(p^-)=(0,0)$,
$$U\big(M^-\cap B_{r_1}^{N^-}(p^-)\big)\equiv\big\{U=(U^1, U^2)\in B_1^k\times B_1^m: \ U^2=0\big\},$$ 
and the Riemannian metric $h^-$ on $B_{r_1}^{N^-}(p^-)$ can be expressed by
$$h^-(U)=\sum_{i,j=1}^{k+m} h_{ij}^-(U) du_i\otimes du_j,\ \forall U\in B_1^k\times B^m_1, $$
and the induced metric of $h^-$ on $M^-\cap B_{r_1}^{N^-}(p^-)$ is given by
$$\displaystyle h^-(U^1,0)=\sum_{i,j=1}^{k} h_{ij}^-(U^1,0)du_i\otimes du_j, \ \forall U^1\in B_1^k.$$

We may assume henceforth that $r_1(N^+)=r_0(M^+)$ in the assumptions ({\bf B}) and ({\bf C}).

\begin{remark} Under the assumptions $({\bf A})$, $({\bf B})$, and $({\bf C})$, it is not hard to see that by choosing a sufficiently small $r_0=r_0(M^+)>0$, under the above local parametrization
of  $B_{r_0}^{N^\pm}(p^\pm)$, the local representations of the Riemannian metrics $h^\pm$ enjoy the following properties:
\begin{equation*}
h^\pm(U)=\sum_{i,j=1}^k h_{ij}^\pm(U^1, U^2) du_i\otimes du_j+\sum_{i,j=k+1}^{k+m} h_{ij}^\pm(U^1, U^2) du_i\otimes du_j, \ \forall U=(U^1,U^2)\in B_1^k\times B_1^m,
\end{equation*}
such that
\begin{equation}\label{asym}
\sum_{i,j=k+1}^{k+m}|h^\pm_{ij}(U^1, U^2)|\le C|U^2|, \forall U=(U^1,U^2)\in B_1^k\times B_1^m,
\end{equation}
for some $C>0$ depending only on $M^\pm$ and $N^\pm$.
\end{remark}

Now we are ready to state a theorem on the local existence of regular solutions to \eqref{HFHM}, whose full proof will be given
in another future work. 
\begin{theorem}\label{existence} Under the assumptions $({\bf A})$, $({\bf B})$, and $({\bf C})$ on $N^\pm$ and $M^\pm$,
for $0<\alpha<1$, let $u_0^\pm\in C^{1+\alpha}(\overline{\Omega^\pm}, N^\pm)$
and $g^\pm=u_0^\pm\big|_{\overline{\Sigma^\pm}}\in C^{1+\alpha}(\overline{\Sigma^\pm}, N^\pm)$ 
be given initial and boundary data such that $u_0^\pm(\Gamma)\subset M^\pm$ satisfies
$u_0^-(x)=\Phi^+(u_0^+(x))$ and $(\frac{\partial u_0^-}{\partial\nu}(x))^T
=D\Phi^+(u^+_0(x))(\frac{\partial u_0^+}{\partial\nu}(x))^T$ for $x\in\Gamma$. Then there exist $T_0>0$, depending on
$\|u_0^\pm\|_{C^{1,\alpha}(\Omega^\pm)}$, and a unique solution
$u^\pm\in C^{1+\alpha,\frac{1+\alpha}2}(Q_{T_0}^\pm, N^\pm)$ of the initial and boundary value problem
\eqref{HFHM}.
\end{theorem}

The proof of Theorem \eqref{existence} is more delicate than the usual proofs of short time smooth solutions to the heat flow of harmonic maps
under the Dirichlet boundary condition (cf. \cite{CL1}, \cite{Hamilton} ) or the free boundary condition (cf. \cite{S3}). It involves to
first show the local existence of regular solutions over small balls, and then patch these local solutions by extending the Schwarz
alternating method on linear parabolic equations to the quasilinear harmonic map heat flows into small neighborhoods
of points in $N^\pm$. For this, we have to overcome major difficulties that arise near the interface $\Gamma$. A detailed proof
will be addressed in a forthcoming work.

In this part, we will indicate a proof of Theorem \ref{existence} when the images of $u^\pm$ is contained
in a single coordinate chart of $N^\pm$. Before doing it, we want to rewrite the system \eqref{HFHM} in an intrinsic form
near a small neighborhood of a point $(x_0,t_0)\in \Gamma_T$ and also derive a generalized energy inequality.

\subsection{Local representation of \eqref{HFHM}}
For $t_0\in (0,T)$ and $x_0\in\Gamma(t_0)$, choose a small $\delta_0>0$, depending on
$\|u^\pm\|_{C^0(Q_T^\pm)}$, such that 
$$u^\pm\big(Q_T^\pm\cap P_{\delta_0}(x_0,t_0)\big)\subset
B_{r_0}^{N^\pm}(p_0^\pm), \ {\rm{with}}\ \ p_0^\pm=u^\pm(x_0,t_0)\in M^\pm.$$
where  $P_{\delta_0}(x_0,t_0)=B_{\delta_0}(x_0)\times (t_0-\delta_0^2, t_0+\delta_0^2)$.
Then, by employing the local representations given by the assumptions
({\bf B}) and ({\bf C}) on $M^\pm, N^\pm$, we can rewrite the 
harmonic heat flow equation \eqref{HFHM} as
\begin{equation}\label{heatflow1}
\begin{cases}
\partial_t U-\Delta U=\Gamma^+(U)(\nabla U, \nabla U) &\ {\rm{in}} \ Q^+_T\cap P_{\delta_0}(x_0,t_0),\\
\partial_t U-\Delta U=\Gamma^-(U)(\nabla U, \nabla U) &\ {\rm{in}} \ Q^{-}_T\cap P_{\delta_0}(x_0,t_0),
\end{cases}
\end{equation}
where 
$U=(U^1,U^2): Q^+_T\cap P_{\delta_0}(x_0,t_0)\mapsto B_1^k\times B_1^m$
is the local representation of $u=u^\pm: Q^+_T\cap P_{\delta_0}(x_0,t_0)\mapsto N$,
and $\Gamma^\pm(\cdot)(\cdot,\cdot)$ is the Christoffel symbol of $N^\pm$. 

Observe that within this  local coordinate system, 
the boundary condition \eqref{HFHM}$_4$
on the free interface $\Gamma_T$ gives rise to
\begin{equation}\label{interface_bdry10}
{U}^2=0,\ {\rm{on}}\ \Gamma_T\cap P_{\delta_0}(x_0,t_0),
\end{equation}
and by \eqref{asym} the boundary condition \eqref{HFHM}$_5$ on the free interface $\Gamma_T$ 
reduces
to 
\begin{equation}\label{interface_bdry20}
\sum_{j=1}^{k}h^+_{ij}(U^1,0)\frac{\partial{(U^1)^j}}{\partial \nu}
=\sum_{j=1}^{k}h_{ij}^-(U^1,0)\frac{\partial{(U^1)^j}}{\partial\nu},  \ 1\le i\le k,
\ {\rm{on}}\ \Gamma_T\cap P_{\delta_0}(x_0,t_0).
\end{equation}

\subsection{Parametrization of domains}

\medskip
\indent Since $\Omega^\pm(t)$ is $t$-dependent over $[0,T]$, in this subsection we will re-parametrize the domains 
and rewrite  \eqref{HFHM} so that it can be viewed as the heat flow of harmonic maps over fixed domain
but with time-dependent metrics on the domain.

Assume that $\Psi(\cdot,t): \Omega\times [0,T]\mapsto\Omega$ is a smooth family of diffeomorphism such that
\begin{equation}\label{diff}
\Psi(x,t)=x,  \forall (x,t)\in\partial\Omega\times [0,T]; 
\ \Psi(\Gamma(t),t)=\Gamma\ {\rm{and}}\ \Psi(\Omega^\pm(t), t)=\Omega^\pm,  \forall\ t\in [0,T].
\end{equation}
For $u^\pm:Q_T^\pm\mapsto N^\pm$, define $\widehat{Q}_T^\pm=\Omega^\pm\times [0,T]$
and $\widehat{u}^\pm:\widehat{Q}_T^\pm\mapsto N^\pm$ through
$${u}^\pm(x,t)=\widehat{u}^\pm(\Psi(x,t), t): \widehat{Q}_T^\pm\mapsto N^\pm.$$

Given that $u^\pm:Q_T^\pm\mapsto N^\pm$ satisfies \eqref{HFHM}, we want to derive the equation
for $\widehat{u}^\pm$ now. To do it,  first set
$$a_{ij}(x,t)=\big(\frac{\partial\Psi_i}{\partial x_\alpha}\frac{\partial\Psi_j}{\partial x_\alpha}\big)(x,t):Q_T
\mapsto \mathbb R^{n\times n},$$
and
$$\widehat{a}_{ij}(y,t)=a_{ij}(x,t): Q_T\mapsto \mathbb R^{n\times n}, \ {\rm{where}}\ (x,t)=\Psi^{-1}(y,t).$$
Then direct calculations imply that
\begin{equation}\label{chain1}
\begin{cases}
\partial_t u^\pm(x,t)=\partial_t \widehat{u}^\pm(\Psi(x,t), t)
+\frac{\partial{\widehat{u}^\pm}}{\partial y_j}(\Psi(x,t), t) \partial_t\Psi_j,\\
\frac{\partial u^\pm}{\partial x_\alpha}=\frac{\partial\widehat{u}^\pm}{\partial y_i}(\Psi(x,t),t)
\frac{\partial\Psi_i}{\partial x_\alpha}, 
\end{cases}
\end{equation}
and
\begin{eqnarray*}
&&\Delta u^\pm(x,t)=\frac{\partial}{\partial x_\alpha}\big(\frac{\partial\widehat{u}^\pm}{\partial y_i}(\Psi(x,t),t)
\frac{\partial\Psi_i}{\partial x_\alpha}\big)\\
&&=\frac{\partial^2\widehat{u}^\pm}{\partial y_i\partial y_j}(\Psi(x,t), t)\frac{\partial\Psi_i}{\partial x_\alpha}
\frac{\partial\Psi_j}{\partial x_\alpha}+\frac{\partial\widehat{u}^\pm}{\partial y_i}(\Psi(x,t),t)\Delta\Psi_i.
\end{eqnarray*}
Hence \eqref{HFHM}$_{1,2}$ becomes
\begin{equation}\label{heatflow2}
\begin{cases}\displaystyle
\partial_t \widehat{u}^+
-\frac{\partial}{\partial y_i}\big(\widehat{a}_{ij} \frac{\partial\widehat{u}^+}{\partial y_j}\big)
= \widehat{a}_{ij}A^+(\widehat{u}^+)\big(\frac{\partial\widehat{u}^+}{\partial y_i},
\frac{\partial\widehat{u}^+}{\partial y_j}\big)
+A_i\frac{\partial\widehat{u}^+}{\partial y_i} & {\rm{in}} \ \widehat{Q}_T^+,\\
\displaystyle\partial_t \widehat{u}^{-}
-\frac{\partial}{\partial y_i}\big(\widehat{a}_{ij} \frac{\partial\widehat{u}^{-}}{\partial y_j}\big)
= \widehat{a}_{ij}A^-(\widehat{u}^-)\big(\frac{\partial\widehat{u}^{-}}{\partial y_i},
\frac{\partial\widehat{u}^{-}}{\partial y_j}\big)
+A_i\frac{\partial\widehat{u}^{-}}{\partial y_i} & {\rm{in}} \ \widehat{Q}_T^-,
\end{cases}
\end{equation}
where
$$A_i(y,t)=\frac{\partial\widehat{a}_{ij}}{\partial y_j}(y,t)
-(\Delta\Psi_i)(\Psi^{-1}(y,t),t)-(\partial_t\Psi_i)(\Psi^{-1}(y,t),t), \ \forall (y,t)\in Q_T.
$$

Observe that 
the boundary condition \eqref{HFHM}$_4$
on the free interface $\Gamma_T$ gives rise to
\begin{equation}\label{interface_bdry1}
\widehat{u}^-(y, t)=\Phi^+(\widehat{u}^+)(y,t), \ \ \forall (y,t)\in \Gamma\times [0,T],
\end{equation}
while the boundary condition \eqref{HFHM}$_5$ on the free interface $\Gamma_T$ gives rise
to 
\begin{equation}\label{interface_bdry2}
\big(\frac{\partial\widehat{u}^-}{\partial \nu}\big)^T(y, t)
=D\Phi^+(\widehat{u}^+)\big(\frac{\partial\widehat{u}^+}{\partial \nu}\big)^T(y, t),  \ \forall (y,t)\in \Gamma\times [0,T],
\end{equation}
where $\nu (=\nu(t))$ is the unit outer normal of $\Gamma$  with respect to the metric 
$\widehat{g}(t)=\widehat{a}_{ij}(y,t)\,dy^idy^j$.


\medskip
First we observe that a sufficiently regular solution of \eqref{HFHM} enjoys a generalized energy inequality. 
For $1<p<\infty$, $T>0$, and an open set $E\subset\mathbb R^n$, denote
$$W^{2,1}_p(E\times [0,T])=\Big\{u\in L^p(E\times [0,T]): \ \partial_t u, \ \nabla^2 u\in L^p(E\times [0,T])\Big\}.$$
We have
\begin{lemma}\label{energy_ineq1} For $T>0$, and $g\in C^{1}({\Sigma^\pm}, N^\pm)$, if $u^\pm\in W^{2,1}_2({Q_T^\pm},  N^\pm)$,
with $\nabla \frac{\partial u^\pm}{\partial t}\in L^2(Q_T^\pm)$,
is a strong solution of \eqref{HFHM}, then there exists constant $C>0$ depending on $\Gamma_T$ such that 
\begin{equation}\label{energy_ineq2}
E(u(t))+\frac14\int_{s}^te^{C(t-\tau)}\big(\int_{\Omega^+(t)} |\partial_t u^+|^2+\int_{\Omega^-(t)} |\partial_t u^-|^2\big)\,dxd\tau\le e^{C(t-s)}E(u(s)), 
\end{equation}
for all $0\le s< t\in [0, T]$.
\end{lemma} 

\pf
Let $\Psi(\cdot, t): \Omega\times [0,T]\mapsto\Omega$ be a smooth family of diffeomorphism given by
\eqref{diff}.
Define $\widehat{u}^\pm: \widehat{Q}^\pm_T\mapsto N^\pm$ by
$$u^\pm(x,t)=\widehat{u}^\pm(\Psi(x,t),t), \ \forall (x,t)\in \widehat{Q}^\pm_T.$$
Then $\widehat{u}^\pm$ solves \eqref{heatflow2} in $\widehat{Q}^\pm_T$, 
\eqref{interface_bdry1} and \eqref{interface_bdry2} on $\Gamma_T$,
and the Dirichlet boundary condition:
\begin{equation}\label{dirichlet}
\widehat{u}^\pm(y,t)=g^\pm(y), \ (y,t)\in\partial\Omega\times [0,T].
\end{equation}
Within this time dependent parametrization, we can write
$$E(u(t))=\frac12\int_{\Omega^+} \widehat{a}_{\alpha\beta}\langle \frac{\partial \widehat{u}^+}{\partial y_\alpha}, \frac{\partial \widehat{u}^+}{\partial y_\beta}\rangle
\,dv_{\widehat{g}}
+\frac12\int_{\Omega^-} \widehat{a}_{\alpha\beta}\langle \frac{\partial \widehat{u}^-}{\partial y_\alpha}, \frac{\partial \widehat{u}^-}{\partial y_\beta}\rangle\,dv_{\widehat{g}},
$$
where $dv_{\widehat{g}}=\sqrt{\widehat{g}}\,dy$, and $\widehat{g}(y,t)={\rm{det}}(\nabla\Psi)(\Psi^{-1}(y,t),t)$.

From $u^\pm\in W^{2,1}_2(Q_T^\pm, N)$ and  $\nabla\frac{\partial u^\pm}{\partial t}\in L^2(Q_T^\pm)$, 
$\widehat{u}^\pm\in W^{2,1}_2(\widehat{Q}^\pm_T, N^\pm)$ and 
$\nabla\frac{\partial \widehat{u}^\pm}{\partial t}\in L^2(\widehat{Q}^\pm_T)$. 
By direct calculations, we have that
\begin{eqnarray*}
\frac{d}{dt}E(u(t))&=&\int_{\Omega^+}\widehat{a}_{\alpha\beta}\langle \frac{\partial}{\partial y_\alpha} (\partial_t{\widehat{u}^+}), 
\frac{\partial \widehat{u}^+}{\partial y_\beta}\rangle\sqrt{\widehat{g}}\,dy+\int_{\Omega^-}\widehat{a}_{\alpha\beta}\langle \frac{\partial}{\partial y_\alpha} (\partial_t{\widehat{u}^-}), 
\frac{\partial \widehat{u}^-}{\partial y_\beta}\rangle\sqrt{\widehat{g}}\,dy\\
&&+\frac12\big(\int_{\Omega^+} \langle \frac{\partial \widehat{u}^+}{\partial y_\alpha}, \frac{\partial \widehat{u}^+}{\partial y_\beta}\rangle
\partial_t(\widehat{a}_{\alpha\beta}\sqrt{{\widehat{g}}})\,dy
+\int_{\Omega^-} \langle \frac{\partial \widehat{u}^-}{\partial y_\alpha}, \frac{\partial \widehat{u}^-}{\partial y_\beta}\rangle
\partial_t(\widehat{a}_{\alpha\beta}\sqrt{{\widehat{g}}})\,dy\big)
\\
&=& I(t)+II(t).
\end{eqnarray*}
It is easy to see that
$$|II(t)|\le CE(u(t)).$$
While, applying the integration by parts, \eqref{heatflow2}, the boundary conditions \eqref{interface_bdry1}, \eqref{interface_bdry2},
and \eqref{dirichlet}, and the
fact that $\partial_t \widehat{u}^-(x,t)=D\Phi^+(u^+)(\partial_t \widehat{u}^{+})(x,t)\in T_{u^\pm(x,t)} M^\pm$ for $(x,t)\in \Gamma_T$,
and $\partial_t \widehat{u}^\pm(x,t)=0$ on $\Sigma^\pm\times [0,T]$, we can show that the boundary contributions on both
$\Gamma$ and $\partial\Omega$ are zeroes. Hence
we can estimate $I$ by
\begin{eqnarray*}
I(t)&=&-\int_{\Omega^+}\langle \partial_t\widehat{u}^+, \frac{\partial}{\partial y_\alpha}(\widehat{a}_{\alpha\beta}\frac{\partial\widehat{u}^+}{\partial y_\beta})\rangle
\,dv_{\widehat{g}}-\int_{\Omega^-}\langle \partial_t\widehat{u}^-, \frac{\partial}{\partial y_\alpha}(\widehat{a}_{\alpha\beta}\frac{\partial\widehat{u}^-}{\partial y_\beta})\rangle
\,dv_{\widehat{g}}\\
&&-
\int_{\Omega^+}\langle \partial_t\widehat{u}^+, \widehat{a}_{\alpha\beta}\frac{\partial\widehat{u}^+}{\partial y_\beta}\rangle \frac{\partial\sqrt{\widehat{g}}}{\partial y_\alpha}
\,dy-\int_{\Omega^-}\langle \partial_t\widehat{u}^-, \widehat{a}_{\alpha\beta}\frac{\partial\widehat{u}^-}{\partial y_\beta}\rangle
\frac{\partial\sqrt{\widehat{g}}}{\partial y_\alpha}\,dy\\
&=&-\big(\int_{\Omega^+}|\partial_t \widehat{u}^+|^2\,dv_{\widehat{g}}+\int_{\Omega^-}|\partial_t \widehat{u}^-|^2\,dv_{\widehat{g}}\big)\\
&&+\big(\int_{\Omega^+}\langle\partial_t \widehat{u}^+, A_i\frac{\partial \widehat{u}^+}{\partial y_i}\rangle\,dv_{\widehat{g}}+\int_{\Omega^-}\langle\partial_t \widehat{u}^-, A_i\frac{\partial \widehat{u}^-}{\partial y_i}\rangle\,dv_{\widehat{g}}\big)\\
&&-\big(\int_{\Omega^+}\langle \partial_t\widehat{u}^+, \widehat{a}_{\alpha\beta}\frac{\partial\widehat{u}^+}{\partial y_\beta}\rangle \frac{\partial\sqrt{\widehat{g}}}{\partial y_\alpha}
\,dy+\int_{\Omega^-}\langle \partial_t\widehat{u}^-, \widehat{a}_{\alpha\beta}\frac{\partial\widehat{u}^-}{\partial y_\beta}\rangle
\frac{\partial\sqrt{\widehat{g}}}{\partial y_\alpha}\,dy\big)\\
&=& III(t)+IV(t)+V(t).
\end{eqnarray*}
It is easy to see that 
$$
|IV(t)|+|V(t)|\le \frac18\big(\int_{\Omega^+}|\partial_t \widehat{u}^+|^2\,dv_{\widehat{g}}+\int_{\Omega^-}|\partial_t \widehat{u}^-|^2\,dv_{\widehat{g}}\big)
+CE(u(t)).
$$
Hence
$$|I(t)|\le -\frac78 \big(\int_{\Omega^+}|\partial_t \widehat{u}^+|^2\,dv_{\widehat{g}}+\int_{\Omega^-}|\partial_t \widehat{u}^-|^2\,dv_{\widehat{g}}\big).
$$
On the other hand, it follows from the chain rule \eqref{chain1} that
\begin{eqnarray*}
&&\int_{\Omega^+}|\partial_t \widehat{u}^+|^2\,dv_{\widehat{g}}+\int_{\Omega^-}|\partial_t \widehat{u}^-|^2\,dv_{\widehat{g}}\\
&&\ge \frac12\big(\int_{\Omega^+(t)}|\partial_t {u}^+|^2\,dx+\int_{\Omega^-(t)}|\partial_t {u}^-|^2\,dx\big)-CE(u(t)).
\end{eqnarray*}
Putting all these estimate together, we obtain
$$
\frac{d}{dt}E(u(t))\le -\frac14 \big(\int_{\Omega^+(t)}|\partial_t {u}^+|^2\,dx+\int_{\Omega^-(t)}|\partial_t {u}^-|^2\,dx\big)+CE(u(t)),
$$
which, combined with Gronwall's inequality, implies \eqref{energy_ineq2}. \qed

\medskip
We will sketch a proof of Theorem \ref{existence} by employing the fixed point argument, under two extra assumptions that
\begin{enumerate}
\item[(i)] the images of 
$u_0^\pm$ is contained in a single coordinate chart, i.e.,
\begin{equation}\label{one-chart}
u_0^\pm(x)\subset B_{r_0}^{N^\pm}(p_0^\pm), \ \forall x\in \overline\Omega, 
\end{equation}
for a pair of points $p_0^\pm\in M^\pm$ that satisfies $p_0^-=\Phi^+(p_0^+)$;
and
\item[(ii)] 
\begin{equation}\label{isom}
\Phi^+:M^+\mapsto M^- \mbox{\ is an isometry}.
\end{equation} 
\end{enumerate}

First we will give some heuristic arguments to indicate that the appropriate function spaces for 
the local existence of regular solutions are
\begin{eqnarray*}
{\mathcal{C}}^{1+\alpha, \frac{1+\alpha}2}_{(u_0,g)}(Q_T, B_{r_0}^{N^\pm}(p_0^\pm))
&=&\Big\{ u:Q_T\mapsto B_{r_0}^{N^\pm}(p_0^\pm):\ u^\pm=u\big|_{Q_T^\pm}\in C^{1+\alpha, \frac{1+\alpha}2}(Q_T^\pm),\\
&& \ \ \ u=u_0 \ {\rm{in}}\ \Omega\times\{0\}, \ \ u=g \ {\rm{on}}\ \partial\Omega\times [0,T],\\ 
&& \ \ \  u^\pm(\Gamma_T)\subset M^\pm,
\ u^-=\Phi^+(u^+),\\
&&\ \ \ (\frac{\partial u^+}{\partial\nu})^T=(D\Phi^+(u^+))^t(\frac{\partial u^-}{\partial\nu})^T \ \ {\rm{on}}\ \ \Gamma_T\Big\},
\end{eqnarray*}
which is equipped with the norm
$$
\big\|u\big\|_{\mathcal{C}^{1+\alpha,\frac{1+\alpha}2}(Q_T)}
=\big\|u^+\big\|_{C^{1+\alpha, \frac{1+\alpha}2}(Q_T^+)}
+\big\|u^-\big\|_{C^{1+\alpha, \frac{1+\alpha}2}(Q_T^-)}.
$$

To see this, assume that $\Gamma(t)\equiv \Gamma$ for $0\le t\le T$. Let
$u^\pm\in {\mathcal{C}}_{(u_0,g)}^{1+\alpha, \frac{1+\alpha}2}(Q_T^\pm, B_{r_0}^{N^\pm}(p_0^\pm))$ be given, and
$U=(U^1, U^2):Q_T\mapsto B_1^k\times B_1^m$ be a local representation of $u^\pm:Q_T^\pm\mapsto B_{r_0}^{N^\pm}(p_0^\pm)$.
Consider $V=(V^1, V^2): Q_T\mapsto B_1^k\times B_1^m$ that is a weak solution of
\begin{equation}\label{heatflow3}
\begin{cases}
\partial_t V-\Delta V=\Gamma^+(U)(\nabla U, \nabla U) &\ {\rm{in}}\ Q_T^+,\\
\partial_t V-\Delta V=\Gamma^-(U)(\nabla U, \nabla U) &\ {\rm{in}} \ Q_T^-,
\end{cases}
\end{equation}
under the initial and boundary condition:
\begin{equation}\label{bdry4}
\begin{cases}
V=U_0  & {\rm{on}}\ \partial_p Q_T,\\
V^2(x^+,t)=V^2(x^-,t)=0,\\
\displaystyle\frac{\partial V^1}{\partial \nu}(x^+,t)
=\frac{\partial V^1}{\partial \nu}(x^-,t), & (x,t)\in \Gamma_T.
\end{cases} 
\end{equation}
Here $U_0:\Omega\mapsto B_1^k\times B_1^m$ is a local representation of $u_0$.

It follows from the regularity of linear parabolic equations
that $V\in C^{1+\alpha, \frac{1+\alpha}2}(Q_T^\pm)$.
Moreover, since
\begin{align*}
\|\Gamma^+(U)(\nabla U, \nabla U)\|_{L^\infty(Q_T^+)}+
\|\Gamma^-(U)(\nabla U, \nabla U)\|_{L^\infty(Q_T^-)}\le C\|\nabla U\|_{L^\infty(Q_T)}^2,
\end{align*}
it follows from the $W^{2,1}_p$-theory of linear parabolic equations that 
$V\in W^{2,1}_p(Q_T^\pm)$ and
\begin{align*}
\big\|V\big\|_{W^{2,1}_p(Q_T^\pm)}\le C(p) \big(\|\nabla U\|_{L^\infty(Q_T)}^2+\|U_0\|_{C^{1+\alpha}(\Omega^\pm)}\big),
\end{align*}
for any $1<p<\infty$. 

By the Sobolev's embedding theorem (\cite{LSU} Lemma II.3.3), we conclude that
$V\in C^{1+\alpha, \frac{1+\alpha}2}(Q_T^\pm)$ and 
\begin{align*} 
\big\|V\big\|_{C^{1+\alpha,\frac{1+\alpha}2}(Q_T^\pm)}\le C(p)\big(\|\nabla U\|_{L^\infty(Q_T)}^2
+\|U_0\|_{C^{1+\alpha}(\Omega^\pm)}\big).
\end{align*}

\medskip
\noindent{\bf Proof of Theorem \ref{existence} under the assumptions \eqref{one-chart} and \eqref{isom}}:

\medskip

For a pair of initial and boundary data $(u_0,g)$ given by Theorem \ref{existence}, 
let $U_0:\partial_p Q_T\mapsto B_1^k\times B_1^m$ be a local representation of $u_0$.
It follows from the assumptions \eqref{one-chart} and \eqref{isom} that 
$u\in \mathcal{C}^{1+\alpha, \frac{1+\alpha}2}_{(u_0,g)}(Q_T, B_{r_0}^{N^\pm}(p_0^\pm))$ if and only if
its local representation $U$ belongs to the space
\begin{eqnarray*}
\mathcal{C}^{1+\alpha, \frac{1+\alpha}2}_{{U_0}}(Q_T, B_1^k\times B_1^m)
&=&\Big\{U=(U^1,U^2)\in C^{1+\alpha, \frac{1+\alpha}2}(Q_T^\pm, B_1^k\times B^m_1): 
U=U_0 \ {\rm{on}}\ \partial_p Q_T,\\
&&\ \ U^2(x^+,t)=U^2(x^-,t)=0, \frac{\partial U^1}{\partial\nu}(x^+,t)=\frac{\partial U^1}{\partial\nu}(x^-,t),
\ (x,t)\in\Gamma_T\Big\}.
\end{eqnarray*}
Now we define $\widehat{U_0}=(\widehat{U_0}^1, \widehat{U_0}^2):Q_T\mapsto B_1^k\times B_1^m$
to the solution of the heat equation in $Q_T^\pm$:
\begin{equation}\label{HE}
\begin{cases}
\partial_t\widehat{U_0}-\Delta \widehat{U_0}=0 & \ {\rm{in}}\ Q_T^\pm,\\
\widehat{U_0}= U_0& \ {\rm{on}}\ \partial_p Q_T,\\
\widehat{U_0}^2(x^+,t)=\widehat{U_0}^2(x^-,t)=0& \ {\rm{on}}\ \Gamma_T,\\
\displaystyle\frac{\partial\widehat{U_0}^1}{\partial\nu}(x^+,t)
=\frac{\partial\widehat{U_0}^1}{\partial\nu}(x^-,t)
&  \ {\rm{on}}\ \Gamma_T.
\end{cases}
\end{equation}
From the condition on ${U_0}$, we know that there exists $\epsilon_0>0$ such that
$$|U_0^1|\le 1-4\epsilon_0, \ |U_0^2|\le 1-4\epsilon_0, \ {\rm{in}}\ \ \Omega.$$
Hence by the maximum principle, we have that
$$|\widehat{U_0}^1|\le 1-2\epsilon_0, \ |\widehat{U_0}^2|\le 1-2\epsilon_0, \ {\rm{in}}\ \ Q_T,$$
and hence $\widehat{U_0}\in \mathcal{C}^{1+\alpha, \frac{1+\alpha}2}_{{U_0}}(Q_T, B_{1-2\epsilon_0}^k\times B_{1-2\epsilon_0}^m)$.

As a consequence, for any $0<\epsilon\le\epsilon_0$, we can see that
$$\mathbb{B}(\widehat{U_0}, \epsilon)=\Big\{U\in \mathcal{C}^{1+\alpha, \frac{1+\alpha}2}_{{U_0}}(Q_T, B_1^k\times B_1^m):
\big\|U-\widehat{U_0}\big\|_{\mathcal{C}^{1+\alpha,\frac{1+\alpha}2}(Q_T)}<\epsilon\Big\}
$$
is a ball in $\mathcal{C}^{1+\alpha, \frac{1+\alpha}2}_{U_0}(Q_T, 
B_1^k\times B_1^m)$ with center $\widehat{U_0}$ and radius $\epsilon$. 

Now we define the solution map $\mathbb{T}: \mathbb{B}(\widehat{U_0},\epsilon)
\mapsto \mathcal{C}^{1+\alpha, \frac{1+\alpha}2}_{\widehat{U_0}}(Q_T, B_1^k\times B_1^m)$ by letting
$V=\mathbb{T}(U), U\in \mathbb{B}(\widehat{U_0},\epsilon)$, be the solution of 
\begin{equation}\label{fix1}
\partial_t V-\Delta V=\begin{cases} \Gamma^+(U)(\nabla U,\nabla U) & {\rm{on}}\ Q_T^+\\
                                                        \Gamma^-(U)(\nabla U,\nabla U) & {\rm{on}}\ Q_T^-,
                                                        \end{cases}
\end{equation}
subject to the initial and boundary condition \eqref{bdry4}.

Now we need the following Lemma.

\begin{lemma} \label{fixpt} There exist $\epsilon>0$ and $T>0$ such that
$\mathbb{T}: \mathbb{B}(\widehat{U_0}, \epsilon)\mapsto \mathbb{B}(\widehat{U_0}, \epsilon)$ is a contractive map,
i.e., for any  $\theta\in (0,1)$, we can find $\epsilon>0$ and $T>0$ such that
\begin{equation}
\big\|\mathbb{T}(U_1)-\mathbb{T}(U_2)\big\|_{\mathcal{C}^{1+\alpha,\frac{1+\alpha}2}(Q_{T})}
\le \theta \big\|U_1-U_2\big\|_{\mathcal{C}^{1+\alpha,\frac{1+\alpha}2}(Q_{T})},
\ \forall U_1, U_2\in \mathbb{B}(\widehat{U_0}, \epsilon).
\end{equation}
Therefore there exists a unique $U\in \mathbb{B}(\widehat{U_0}, \epsilon)$ such that $U=\mathbb{T}(U)$.
In particular, if $u^\pm=u\big|_{Q_{T}^\pm}: Q_{T}^\pm\mapsto N^\pm$ has $U$ as its
local representation, then $u$ is a unique regular solution of \eqref{HFHM} in $Q_{T}$.
\end{lemma} 

\pf
For $\mathbb{U}\in \mathbb{B}(\widehat{U_0},\epsilon)$, since $V-\widehat{U_0}$ satisfies
\begin{equation}\label{fix2}
\partial_t (V-\widehat{U_0})-\Delta (V-\widehat{U_0})=\begin{cases} \Gamma^+(U)(\nabla U,\nabla U) & {\rm{on}}\ Q_T^+\\
                                                        \Gamma^-(U)(\nabla U,\nabla U) & {\rm{on}}\ Q_T^-,
                                                        \end{cases}
\end{equation}
and 
$$\begin{cases}
V-\widehat{U_0}=0 &\ {\rm{on}}\ \partial_p Q_T,\\
(V-\widehat{U_0})^2(x^+,t)=(V-\widehat{U_0})^2(x^-,t)=0 &\ {\rm{on}}\ \Gamma_T,\\
\frac{\partial(V-\widehat{U_0})^1}{\partial\nu}(x^+,t)
=\frac{\partial(V-\widehat{U_0})^1}{\partial\nu}
(x^-,t) &\ {\rm{on}}\ \Gamma_T.
\end{cases}
$$
Hence, similar to the earlier discussion, we have that for some $p=p(\alpha)>n+2$, 
\begin{eqnarray*}
\big\|V-\widehat{U_0}\big\|_{\mathcal{C}^{1+\alpha, \frac{1+\alpha}2}(Q_T)}
&\le& C\big\|V-\widehat{U_0}\big\|_{W^{2,1}_p(Q_T^\pm)}\\
&\le& C\big\||\nabla U|^2\big\|_{L^p(Q_T)}\\
&\le& C\big\|\nabla U\big\|_{L^\infty(Q_T)}^2 T^{\frac{1}{p}}\\
&\le & C(\|\widehat{U_0}\|_{\mathcal{C}^{1+\alpha, \frac{1+\alpha}2}(Q_T)}+\epsilon)^2 T^{\frac{1}{p}}\\
&< & \epsilon,
\end{eqnarray*} 
provided we choose a sufficiently small $T=T_0>0$, depending only on $U_0$ and $\alpha$.
Hence $V=\mathbb{T}({U})\in \mathbb{B}(\widehat{U_0},\epsilon)$. 

For $i=1, 2$, let  $U_i\in \mathbb{B}(\widehat{U_0},\epsilon)$ and $V_i=\mathbb{T}(U_i)$. Then
\begin{equation}\label{fix3}
\partial_t (U_1-U_2)-\Delta (U_1-U_2)=\begin{cases} \Gamma^+(U_1)(\nabla U_1,\nabla U_1)
-\Gamma^+(U_2)(\nabla U_2,\nabla U_2) & {\rm{on}}\ Q_T^+\\
\Gamma^-(U_1)(\nabla U_1,\nabla U_1)-\Gamma^-(U_2)(\nabla U_2,\nabla U_2) & {\rm{on}}\ Q_T^-,
\end{cases}
\end{equation}
and 
$$\begin{cases}
U_1-{U_2}=0 &\ {\rm{on}}\ \partial_p Q_T,\\
(U_1-{U_2})^2(x^+,t)=(U_1-{U_2})^2(x^-,t)=0 &\ {\rm{on}}\ \Gamma_T,\\
\frac{\partial(U_1-{U_2})^1}{\partial\nu}(x^+,t)
=\frac{\partial(U_1-{U_2})^1}{\partial\nu}
(x^-,t) &\ {\rm{on}}\ \Gamma_T.
\end{cases}
$$
Hence we can conclude  that  for any $\theta\in (0,1)$ such that for $p=p(\alpha)>n+2$,
 \begin{eqnarray*}
&&\big\|U_1-{U_2}\big\|_{\mathcal{C}^{1+\alpha, \frac{1+\alpha}2}(Q_T)}\\
&\le& C\big\|U_1-{U_2}\big\|_{W^{2,1}_p(Q_T^\pm)}\\
&\le& C\big\|(|\nabla U_1|+|\nabla U_2|)^2(U_1-U_2)+(|\nabla U_1|+|\nabla U_2|)|\nabla(U_1-U_2)|\big\|_{L^p(Q_T)}\\
&\le& C\big\||\nabla U_1|+|\nabla U_2|\big\|_{L^p(Q_T)}^2 \|U_1-U_2\|_{L^\infty(Q_T)}
+C\big\||\nabla U_1|+|\nabla U_2|\big\|_{L^p(Q_T)}\|\nabla(U_1-U_2)\|_{L^\infty(Q_T)}\\
&\le & C(1+\|\widehat{U_0}\|_{\mathcal{C}^{1+\alpha, \frac{1+\alpha}2}(Q_T)}^2) T^{\frac{1}{p}}
\big\|U_1-{U_2}\big\|_{\mathcal{C}^{1+\alpha, \frac{1+\alpha}2}(Q_T)}\\
&< & \theta \big\|U_1-{U_2}\big\|_{\mathcal{C}^{1+\alpha, \frac{1+\alpha}2}(Q_T)},
\end{eqnarray*} 
provided $T=T_0>0$ is chosen so that 
$$C(1+\|\widehat{U_0}\|_{\mathcal{C}^{1+\alpha, \frac{1+\alpha}2}(Q_T)}^2) T_0^{\frac{1}{p}}\le \theta.$$
This completes the proof of both Lemma \ref{fixpt} and Theorem \ref{existence} under the assumptions
\eqref{one-chart} and \eqref{isom}.
\qed

\section{Boundary monotonicity inequality of \eqref{HFHM}}
\setcounter{equation}{0}
\setcounter{theorem}{0}

In this section, we will derive a boundary monotonicity
inequality on \eqref{HFHM}, analogous to Struwe's monotonicity formula, which may have its own interest.

To simplify the presentation, we assume that 
$$\Omega=\mathbb R^n, \ T>0, \ {\rm{and}}\ \Gamma(t)=\Gamma=\partial\mathbb R^n_+ \  {\rm{for}}\ 0\le t\le T.$$

Let $u^\pm:\mathbb R^n_\pm\times [0,+\infty)\to N^\pm$, with $u^\pm(x, t)\in M^\pm$ for $(x,t)\in \partial \mathbb R^n_+\times (0,\infty)$, satisfy
\begin{equation}\label{heatflow5}
\begin{cases}
\partial_t u^+-\Delta u^+=A^+(u^+)(\nabla u^+, \nabla u^+)  & {\rm{in}}\ \mathbb R^n_+\times (0,+\infty)\\
\partial_t u^- -\Delta u^-=A^-(u^-)(\nabla u^-, \nabla u^-)  & {\rm{in}}\ \mathbb R^n_-\times (0,+\infty)\\
\Phi^+(u^+)=u^- & {\rm{in}}\ \partial \mathbb R^n_+\times (0,+\infty)\\
(\frac{\partial u^+}{\partial x_n})^T=\Big(D\Phi^+(u^+)\Big)^t\Big[(\frac{\partial u^-}{\partial x_n})^T\Big]
 & {\rm{on}}\ \partial\mathbb R^n_+\times (0,+\infty).
\end{cases}
\end{equation}

For $(x_0,t_0)\in\mathbb R^n\times (0,+\infty)$ and $0<R\le\sqrt{t_0}$,  let
$$G_{(x_0,t_0)}(x,t)=\frac{1}{(4\pi(t_0-t))^{\frac{n}2}}e^{-\frac{|x-x_0|^2}{4(t_0-t)}},
\ (x,t)\in\mathbb R^n\times (0,t_0)
$$
denoet the backward heat kernel on $\mathbb R^n$.
Set
$$E(u^\pm; (x_0,t_0), R)=R^2\Big[\int_{\mathbb R^n_+\times \{t_0-R^2\}}|\nabla u^+|^2G_{(x_0,t_0)}(x,t)\,dx
+\int_{\mathbb R^n_{-}\times \{t_0-R^2\}}|\nabla u^-|^2G_{(x_0,t_0)}(x,t)\,dx\Big].
$$
\begin{lemma} \label{bdry_mono}
Suppose that $(x_0,t_0)=(0,0)\in\partial\mathbb{R}^n_+\times (-\infty, 0]$ and $u^\pm \in C^2(\overline{\mathbb R^n_\pm}\times (-\infty, 0], N^\pm)$ is a solution to the system (\ref{heatflow5}). Then
\begin{equation}
\frac{d}{dR}E(R)\ge 0. \label{bdry_mono1}
\end{equation}
\end{lemma}
\pf Write $G(x,t)$ for $G_{(0,0)}(x,t)$ and define 
$$u_R^\pm(x,t)=u^\pm(Rx,R^2 t), \ (x,t)\in\mathbb R^n\times (-\infty,0].$$
It is easy to see that
$$E(u^\pm; R)=E(u^\pm_R; 1).$$
For simplicity, we only verify (\ref{bdry_mono1}) at $R=1$. Since
$$\frac{d}{dR}\big|_{R=1} u_R^\pm=x\cdot\nabla u^\pm-2\partial_t u^\pm,$$
we have
\begin{eqnarray*}
&&\frac{d}{dR}\Big|_{R=1}E(u^\pm; R)=\frac{d}{dR}\Big|_{R=1}E(u^\pm_R; 1)\\
&=&2\int_{\mathbb R^n_+\times \{-1\}} \nabla u^+\cdot\nabla(x\cdot\nabla u^+-2\partial_t u^+)
e^{-\frac{|x|^2}4}\\
&&+2\int_{\mathbb R^n_{-}\times \{-1\}} \nabla u^-\cdot\nabla(x\cdot\nabla u^--2\partial_t u^-)
e^{-\frac{|x|^2}4}\\
&=&-2\Big[\int_{\mathbb R^n_+\times \{-1\}} \nabla\cdot(\nabla u^+ e^{-\frac{|x|^2}4})\cdot (x\cdot\nabla u^+-2\partial_t u^+)\\
&&+\int_{\mathbb R^n_{-}\times \{-1\}} \nabla\cdot(\nabla u^-e^{-\frac{|x|^2}4})\cdot (x\cdot\nabla u^--2\partial_t u^-)\Big]\\
&&-2\Big[\int_{\partial \mathbb R^n_+\times \{-1\}} \frac{\partial u^+}{\partial x_n} 
\cdot (x\cdot\nabla u^+-2\partial_t u^+)e^{-\frac{|x|^2}4}\\
&&-\int_{\partial \mathbb R^n_{-}\times \{-1\}} \frac{\partial u^-}{\partial x_n}
\cdot (x\cdot\nabla u^--2\partial_t u^-)e^{-\frac{|x|^2}4}\Big].\\
\end{eqnarray*}
Since 
$$\nabla\cdot(\nabla u^\pm e^{-\frac{|x|^2}4})
=\Delta u^\pm-\frac{1}{2}x\cdot\nabla u^\pm
=\partial_t u^\pm-A^\pm(u^\pm)(\nabla u^\pm, u^\pm)-\frac{1}{2}x\cdot\nabla u^\pm
$$
and
$$A^\pm(u^\pm)(\nabla u^\pm, u^\pm)\cdot 
(x\cdot\nabla u^\pm-2\partial_t u^\pm)=0,$$
we have
\begin{eqnarray*}
&&-2\Big[\int_{\mathbb R^n_+\times \{-1\}} \nabla\cdot(\nabla u^+ e^{-\frac{|x|^2}4})\cdot (x\cdot\nabla u^+-2\partial_t u^+)\\
&&+\int_{\mathbb R^n_{-}\times \{-1\}} \nabla\cdot(\nabla u^-e^{-\frac{|x|^2}4})\cdot (x\cdot\nabla u^--2\partial_t u^-)\Big]\\
&=&\big[\int_{\mathbb R^n_+\times \{-1\}} |x\cdot\nabla u^+-2\partial_t u^+|^2
e^{-\frac{|x|^2}4}+\int_{\mathbb R^n_{-}\times \{-1\}} |x\cdot\nabla u^--2\partial_t u^-|^2
 e^{-\frac{|x|^2}4}\big].
\end{eqnarray*}
Since $x=(x',0)$ for $x\in\partial \mathbb R^n_\pm$, and 
$$u^\pm\big(\partial\mathbb R^n_\pm\times (-\infty,0)\big)\subset M^\pm,$$
we have 
$$x\cdot\nabla u^\pm-2\partial_t u^\pm\big|_{\partial\mathbb R^n_\pm\times (-\infty, 0)}
\in T_{u^\pm(x,t)}M^\pm$$
so that
$$\frac{\partial u^\pm}{\partial x_n}
\cdot (x\cdot\nabla u^\pm-2\partial_t u^\pm)
=\Big(\frac{\partial u^\pm}{\partial x_n}\Big)^T
\cdot (x\cdot\nabla u^\pm-2\partial_t u^\pm)
\ \ {\rm{on}}\ \ {\partial\mathbb R^n_\pm\times (-\infty, 0)}.
$$
Since
$$u^-(x,t)=\Phi^+(u^+(x,t)) \ \ {\rm{on}}\ \ {\partial\mathbb R^n_\pm\times (-\infty, 0)},$$
we have
$$\nabla_{\rm{tan}} u^-(x,t)=D\Phi^+(u^+(x,t))\nabla_{\rm{tan}} u^+(x,t),
\ \partial_t u^-(x,t)=D\Phi^+(u^+(x,t))\partial_t u^+(x,t)$$
$ {\rm{on}}\ \ {\partial\mathbb R^n_\pm\times (-\infty, 0)}$ and hence
$$x\cdot\nabla u^-(x,t)-2\partial_t u^-(x,t)
=D\Phi^+(u^+(x,t))\big(x\cdot\nabla u^+(x,t)-2\partial_t u^+(x,t)\big)
\ {\rm{on}}\ \ {\partial\mathbb R^n_\pm\times (-\infty, 0)}.$$
Therefore we have
\begin{eqnarray*}
&&\big[\int_{\partial \mathbb R^n_+\times \{-1\}} \frac{\partial u^+}{\partial x_n} 
\cdot (x\cdot\nabla u^+-2\partial_t u^+)e^{-\frac{|x|^2}4}\\
&&-\int_{\partial \mathbb R^n_{-}\times \{-1\}} \frac{\partial u^-}{\partial x_n}
\cdot (x\cdot\nabla u^--2\partial_t u^-)e^{-\frac{|x|^2}4}\Big]\\
&=&\int_{\partial\mathbb R^n_+\times \{-1\}}
\big[(\frac{\partial u^+}{\partial x_n})^T-(D\Phi^+(u^+(x,t)))^t
(\frac{\partial u^-}{\partial x_n})^T\big]
\cdot (x\cdot\nabla u^+-2\partial_t u^+)e^{-\frac{|x|^2}4}=0,
\end{eqnarray*}
where we have used the boundary condition (\ref{heatflow5})$_5$ in the last step. Putting all
these calculations together, we obtain
\begin{eqnarray}\label{mono_id2}
&&\frac{d}{dR}\Big|_{R=1}E(u^\pm; R)\nonumber\\
&=&\big[\int_{\mathbb R^n_+\times \{-1\}} |x\cdot\nabla u^+-2\partial_t u^+|^2
e^{-\frac{|x|^2}4}+\int_{\mathbb R^n_{-1}\times \{-1\}} |x\cdot\nabla u^--2\partial_t u^-|^2
 e^{-\frac{|x|^2}4}\big]\ge 0.
\end{eqnarray}
This completes the proof. \qed

\end{document}